\documentclass{article}

\usepackage{arxiv}
\usepackage[utf8]{inputenc} 
\usepackage[T1]{fontenc}    
\usepackage{hyperref}       
\usepackage{url}            
\usepackage{booktabs}       
\usepackage{amsfonts}       
\usepackage{nicefrac}       
\usepackage{microtype}      
\usepackage{tikz}
\usepackage{amssymb}
\usepackage{amsmath}
\usepackage{tcolorbox}
\usepackage{fullpage} 

\tcbuselibrary{listings,skins,theorems}

\lstdefinelanguage{Maple}%
{aboveskip={0pt},belowskip={0pt},basicstyle=\bfseries,
morekeywords={and,assuming,break,by,catch,description,do,done,%
elif,else,end,error,export,fi,finally,for,from,global,if,%
implies,in,intersect,local,minus,mod,module,next,not,od,%
option,options,or,proc,quit,read,return,save,stop,subset,then,%
to,try,union,use,uses,while,xor},%
sensitive=true,%
morecomment=[l]\#,%
morestring=[b]",%
morestring=[d]"%
}[keywords,comments,strings]%

\makeatletter
\tcbset{every listing line/.code={%
  \def\lst@NewLine{%
      \ifx\lst@OutputBox\@gobble\else
          \par\noindent \hbox{}#1%
      \fi
      \global\advance\lst@newlines\m@ne
      \lst@newlinetrue}}}
\makeatother

\newtcblisting{code}[1]{
           colback=white,
           skin=bicolor,
           colbacklower=white,
           colupper=red!40!black,
           collower=blue,
           frame style={draw=white,left color=white,right color=white},
           fontupper=\bfseries,
           fontlower=\bfseries,
           width=\linewidth,
           boxrule=0mm,
           outer arc=1mm,arc=1mm,
           leftupper=0cm,leftlower=0cm,rightupper=0cm,rightlower=0cm,
           top=0mm,bottom=0mm,middle=0mm,
           center lower,
           nobeforeafter, 
           listing and comment,
           ams nodisplayskip lower,
           every listing line={\textcolor{red!40!black}{\ttfamily> }}
           }

\newtheorem{prop}{Proposition}[section]
\newtheorem{problem}[prop]{Problem}
\newtheorem{thm}[prop]{Theorem}

\title{Some arithmetical problems that are obtained by analyzing proofs and infinite graphs}

\author{
  Lorenzo Sauras-Altuzarra \\
  Research Unit of Computational Logic \\
  TU Wien, Austria \\
  \texttt{lorenzo.sauras@tuwien.ac.at} \\
}

\begin{document}
\maketitle

\begin{abstract}

Applying Baaz's Generalization Method and a new technique to, respectively, proofs and denumerable simple graphs, diverse arithmetical patterns are observed. In particular, sufficient conditions for a number to be a divisor of a Fermat number are provided. The accuracy of such observations is asked in several subsequent problems.

\end{abstract}

\keywords{Cayley graph \and Collatz's Problem \and Fermat number \and Graph sequence \and Integer sequence \and Open problem \and Proof analysis \and Proof complexity}

\section{Introduction}

Roughly speaking, Baaz's Generalization Method (see \cite{Baaz}) is a procedure that, given a proof of a certain theorem, allows to obtain a more general proof (and, in particular, a more general theorem). In this paper the idea is applied to six proofs that fifth Fermat number is divisible by 641, deriving in several sufficient conditions for a number to be a factor of a Fermat number, remarkably:
\begin{itemize}
    \item $ k \cdot 2^s + 1 \ | \ F_n $, for every $ k , n , r , s \in \mathbb{N}^+ $ such that $ r \cdot s \leq 2^{n - 1} $ and $ k \cdot 2^s + 1 \ | \ k^{2 \cdot r} + 2^{2^n - 2 \cdot r \cdot s} $ (Theorem \ref{Particularizacion}),
    \item $ 2^{2^n - 4 \cdot ( n + 2 )} + i^4 \ | \ F_n $, for every $ i , n \in \mathbb{N}^+ $ such that $ n > 4 $ and $ i \cdot 2^{n + 2} + 1 = 2^{2^n - 4 \cdot ( n + 2 )} + i^4 $ (Theorem \ref{Particularizacion3}),
    \item $ i \ | \ F_n $, for every $ c , i , n \in \mathbb{N}^+ $ such that $ i \ | \ ( 2^{2^{n - 1}} - i \cdot c )^2 + 1 $ (Theorem \ref{Particularizacion4}).
\end{itemize}

In addition, a bijection between the so-called vertex-by-vertex increasing graph sequences and certain integer sequences is provided, together with various examples of its application.

Throughout the text, numerous problems are posed. In many cases such problems are accompanied by Maple implementations that can help the reader to experiment with the described concepts.

This article is written with the intention of giving an idea about how both techniques can be a useful tool for pattern recognition in diverse mathematical areas.

\section{Generalizations of proofs}

In this section Baaz's Generalization Method is explained and applied to six proofs that fifth Fermat number is divisible by 641, obtaining some interesting generalizations and questions.

\newpage

\subsection{Baaz's Generalization Method}{\label{SectionExplain}}

In mathematics, examples are very important, but not all of them are equally good. Intuitively, when accompanying a theorem, the more it reflects the potential of such result, the better is an example. In fact, if an example sufficiently represents the essence of a theorem, then it can be almost so instructive as the proof itself. Mathematical teaching via examples have diverse endorsements, remarkably Babylonian mathematics, which mainly consisted of collections of examples (see Chapter 3 of \cite{VanderWaerden}).

Baaz's Generalization Method formalizes a way of measuring the "goodness" of examples. Indeed, from a concrete example $ E $ of certain universal theorem $ T $ (i.e., a theorem that can be expressed by a universal formula, see \cite{Weisstein2}), it generates another universal theorem $ t ( E ) $, with its corresponding proof. A subsequent comparison between $ T $ and $ t ( E ) $ may show how approximate was $ E $ to $ T $.

But there is yet another possible use of this procedure: if applied to an answer to a concrete case of an open problem (as a proof that 641 divides the fifth Fermat number), it will output a result that can be particularized to a partial answer of such question (as a sufficient condition for a number to be a divisor of an arbitrary Fermat number).

An explanation for this algorithm in full generality can be found in \cite{Baaz}. For elementary number theory, which is the case at hand in this article, it works as follows.

\medskip

\noindent\begin{tabular}{|p{7.5cm}|p{7.5cm}|}
\hline
\textbf{Baaz's Generalization Method} & Example \\
\hline
Input: proof of a universal formula (i.e., tree of predicates that are connected by implications). & Input:

$ 641 \ | \ \underbrace{5^4 + 2^4}_{641} \ \Rightarrow \ 641 \ | \ 5^4 \cdot 2^{28} + 2^{32} \ \Rightarrow $

$ 641 \ | \ ( 5^4 \cdot 2^{28} - 1 ) + ( 2^{32} + 1 ) $. \hfill $ [ 1 ] $

$ \underbrace{5 \cdot 2^7 + 1}_{641} \ | \ 5^4 \cdot 2^{28} - 1 \ \overset{[ 1 ]}{\Rightarrow} \ 641 \ | \ 2^{32} + 1 $. \hfill $ \square $ \\
\hline
1. For every leaf (i.e., vertex without predecessors), replace every constant by a variable (without repeating them) (there is no need to keep the operations, but keep the rest of the relations). & 1.

$ 641 \ | \ 5^4 + 2^4 \ \mapsto \ a_0 \ | \ b_0 $,

$ 5 \cdot 2^7 + 1 \ | \ 5^4 \cdot 2^{28} - 1 \ \mapsto \ c_0 \ | \ d_0 $. \\
\hline
2. For every implication, replace every constant by a variable (without repeating them) (keep the operations that are necessary to justify the step and the rest of the relations). & 2.

$ 1^\textrm{st} \ \Rightarrow $: $ [ a_1 \ | \ b_1 + c_1 \ \Rightarrow \ a_1 \ | \ b_1 \cdot d_1 + c_1 \cdot d_1 ] $,

$ 2^\textrm{nd} \ \Rightarrow $: $ [ a_2 \ | \ b_2 + c_2 \ \Rightarrow \ a_2 \ | \ ( b_2 - d_2 ) + ( c_2 + d_2 ) ] $,

$ 3^\textrm{rd} \ \Rightarrow $: $ [ [ a_3 \ | \ b_3 + c_3 \ \wedge \ a_3 \ | \ b_3 ] \ \Rightarrow \ a_3 \ | \ c_3 ] $. \\
\hline
3. Minimize the number of variables, by simultaneously unifying all pairs of predicates that are assigned to the same vertex (all the relations (and, in particular, all the operations) must be kept). & 3.

$ \{ a_0 \ | \ b_0 , a_1 \ | \ b_1 + c_1 \} \ \mapsto \ A \ | \ D + B $,

$ \{ a_1 \ | \ b_1 \cdot d_1 + c_1 \cdot d_1 , a_2 \ | \ b_2 + c_2 \} \ \mapsto \ A \ | \ D \cdot C + B \cdot C $,

$ \{ a_2 \ | \ ( b_2 - d_2 ) + ( c_2 + d_2 ) , a_3 \ | \ b_3 + c_3 \} \ \mapsto $

$ A \ | \ ( D \cdot C - E ) + ( B \cdot C + E ) $,

$ \{ c_0 \ | \ d_0 , a_3 \ | \ b_3 \} \ \mapsto \ A \ | \ D \cdot C - E $. \\
\hline
Output: generalized proof (and, in particular, generalized theorem, whose hypotheses are the generalized leaves and whose thesis is the generalized root). & Output:

$ A \ | \ D + B \ \Rightarrow \ A \ | \ D \cdot C + B \cdot C \ \Rightarrow $

$ A \ | \ ( D \cdot C - E ) + ( B \cdot C + E ) $. \hfill $ [ 1 ] $

$ A \ | \ D \cdot C - E \ \overset{[ 1 ]}{\Rightarrow} \ A \ | \ B \cdot C + E $ \hfill $ \square $

(generalized theorem: if $ A \ | \ D + B $ and $ A \ | \ D \cdot C - E $, then $ A \ | \ B \cdot C + E $). \\
\hline
\end{tabular}

\subsection{Proofs that the fifth Fermat number is divisible by 641 and their generalizations}

Recall that, given $ r \in \mathbb{N} $, $ r $ is a \textbf{Fermat number} if, and only if, $ r = 2^{2^n} + 1 $, for some $ n \in \mathbb{N} $. From now on, the strictly increasing integer sequence whose image is the set of Fermat numbers will be denoted by $ F $.

At the present moment, one of the most important problems in arithmetic is the following one.

\begin{problem}
    Does there exist $ n \in \mathbb{N} $ such that $ n > 4 $ and $ F_n $ is a prime number?
\end{problem}

It will be convenient to keep in mind the following results.

\begin{prop}
    $ F_0 , ... , F_4 $ are prime numbers (see \cite{Weisstein}).
\end{prop}

\begin{thm}[Euler \& Lucas]
    If $ n \in \mathbb{N} \backslash \{ 0 , 1 \} $, then $ k \cdot 2^{n + 2} + 1 \ | \ F_n $, for some $ k \in \mathbb{N}^+ $ (see \cite{Weisstein}).
\end{thm}

\begin{thm}{\label{RecDef}}
    $ F_0 = 3 $ and $ F_{n + 1} = \displaystyle\prod_{k = 0}^n ( F_k ) + 2 $, for every $ n \in \mathbb{N} $ (see \cite{Weisstein}).
\end{thm}

\subsubsection{The six basic proofs}

\begin{thm}
    $ 641 \ | \ F_5 $.
\end{thm}

\underline{Proof I} (Bennet \& Kraïtchik, see page 19 of \cite{Krizek})

$ 641 \ | \ \underbrace{5^4 + 2^4}_{641} \ \Rightarrow \ 641 \ | \ 5^4 \cdot 2^{28} + 2^{32} \ \Rightarrow \ 641 \ | \ ( 5^4 \cdot 2^{28} - 1 ) + \underbrace{( 2^{32} + 1 )}_{F_5} $. \hfill $ [ 1 ] $

$ \underbrace{5 \cdot 2^7 + 1}_{641} \ | \ 5^4 \cdot 2^{28} - 1 \ \overset{[ 1 ]}{\Rightarrow} \ 641 \ | \ F_5 $. \hfill $ \square $

\medskip

\underline{Proof II} (see page 160 of \cite{Coppel})

$ \underbrace{5 \cdot 2^7 + 1}_{641} \equiv 0 \ ( \textrm{mod}. \ 641 ) \ \Rightarrow \ \underbrace{5 \cdot 2^7 + 1 - 2^4}_{5^4} \equiv - 2^4 \ ( \textrm{mod}. \ 641 ) \ \Rightarrow \ ( \underbrace{5 \cdot 2^7}_{640} )^4 \equiv - 2^{32} \ ( \textrm{mod}. \ 641 ) $. \hfill $ [ 1 ] $

$ 640 \equiv - 1 \ ( \textrm{mod}. \ 641 ) \ \Rightarrow \ 640^4 \equiv 1 \ ( \textrm{mod}. \ 641 ) \ \overset{[ 1 ]}{\Rightarrow} \ - 2^{32} \equiv 1 \ ( \textrm{mod}. \ 641 ) \ \Rightarrow \ 641 \ | \ \underbrace{2^{32} + 1}_{F_5} $. \hfill $ \square $

\medskip

\underline{Proof III} (Kraïtchik, see page 39 of \cite{Krizek})

$ 4 \ \textrm{is even} \ \Rightarrow \ \displaystyle\sum_{k = 0}^4 \left( \displaystyle\binom{4}{k} \cdot 641^{4 - k} \cdot ( - 1 )^k \right) = \displaystyle\sum_{k = 0}^{4 - 1} \left( \displaystyle\binom{4}{k} \cdot 641^{4 - k} \cdot ( - 1 )^k \right) + 1 $. \hfill $ [ 1 ] $

$ 5^4 = 641 - ( 641 - 5^4 ) \ \Rightarrow \ 5^4 \cdot 2^{28} = 641 \cdot 2^{28} - ( 641 - 5^4 ) \cdot 2^{28} $. \hfill $ [ 2 ] $
    
$ 5 \cdot 2^7 = 641 + ( - 1 ) \ \Rightarrow \ 5^4 \cdot 2^{28} = ( 641 + ( - 1 ) )^4 \ \overset{\tiny \begin{array}{c} \textrm{Binomial} \\ \textrm{theorem} \end{array}}{\Rightarrow} $
    
$ 5^4 \cdot 2^{28} = \displaystyle\sum_{k = 0}^4 \left( \displaystyle\binom{4}{k} \cdot 641^{4 - k} \cdot ( - 1 )^k \right) \ \overset{[ 1 ]}{\Rightarrow} $

$ 5^4 \cdot 2^{28} = \displaystyle\sum_{k = 0}^{4 - 1} \left( \displaystyle\binom{4}{k} \cdot 641^{4 - k} \cdot ( - 1 )^k \right) + 1 \ \Rightarrow $
    
$ 5^4 \cdot 2^{28} = 641 \cdot \displaystyle\sum_{k = 0}^{4 - 1} \left( \displaystyle\binom{4}{k} \cdot 641^{4 - k - 1} \cdot ( - 1 )^k \right) + 1 \ \overset{[ 2 ]}{\Rightarrow} $
    
$ 641 \cdot 2^{28} - ( 641 - 5^4 ) \cdot 2^{28} = 641 \cdot \displaystyle\sum_{k = 0}^3 \left( \displaystyle\binom{4}{k} \cdot 641^{4 - k - 1} \cdot ( - 1 )^k \right) + 1 \ \Rightarrow $
    
$ 641 \cdot \left( 2^{28} - \displaystyle\sum_{k = 0}^3 \left( \displaystyle\binom{4}{k} \cdot 641^{4 - k - 1} \cdot ( - 1 )^k \right) \right) = \underbrace{( \underbrace{641}_{5^4 + 2^4} - 5^4 ) \cdot 2^{28} + 1}_{F_5} \ \Rightarrow $

$ 641 \ | \ F_5 $. \hfill $ \square $

\medskip

\underline{Proof IV} (Bennet, see Proposition 4.17 of \cite{Krizek})

$ 2^7 - 5^3 - 3 = 0 \ \Rightarrow \ - 5 \cdot 3 + 5 \cdot ( 2^7 - 5^3 ) = 0 \ \Rightarrow \ 1 + 2^7 \cdot 5 - 5^4 = 2^4 $. \hfill $ [ 1 ] $

$ 4 + 28 = 4 + 28 \ \Rightarrow \ 2^{4 + 28} = 2^{4 + 28} \ \Rightarrow \ 2^{4 + 28} + 1 = 2^{4 + 28} + 1 \ \overset{[ 1 ]}{\Rightarrow} $

$ \underbrace{2^{4 + 28} + 1}_{F_5} = ( \underbrace{1 + 2^7 \cdot 5}_{641} - 5^4 ) \cdot 2^{28} + 1 \ \Rightarrow $

$ F_5 = 641 \cdot 2^{28} + ( 1 - 5^4 \cdot 2^{28} ) \ \Rightarrow $

$ F_5 = 641 \cdot 2^{28} + ( 1 + 5^2 \cdot 2^{14} ) \cdot \underbrace{( 1 + 5 \cdot 2^7 )}_{641} \cdot ( 1 - 5 \cdot 2^7 ) \ \Rightarrow $

$ F_5 = 641 \cdot ( 2^{28} + ( 1 + 5^2 \cdot 2^{14} ) \cdot ( 1 - 5 \cdot 2^7 ) ) \ \Rightarrow $

$ 641 \ | \ F_5 $. \hfill $ \square $

\medskip

\underline{Proof V} (Sleziak, see \underline{\url{https://math.stackexchange.com/q/149658}})

$ 641 \ | \ 641 \ \Rightarrow \ 641 \ | \ 641 \cdot 37 \ \Rightarrow \ 641 \ | \ 154^2 - \underbrace{( 154^2 - 641 \cdot 37 )}_{- 1} \ \Rightarrow $

$ 641 \ | \ ( 154^2 - 2^{32} ) + \underbrace{( 1 + 2^{32} )}_{F_5} $. \hfill $ [ 1 ] $

$ 641 \ | \ 640 - ( - 1 ) \ \Rightarrow \ 641 \ | \ 640 \cdot 102 - ( - 1 ) \cdot 102 \ \Rightarrow \ 641 \ | \ \underbrace{( 256 + 640 \cdot 102 )}_{2^{16}} - \underbrace{( 256 + ( - 1 ) \cdot 102 )}_{154} \ \Rightarrow $

$ 641 \ | \ 2^{32} - 154^2 \ \Rightarrow \ 641 \ | \ 154^2 - 2^{32} \ \overset{[ 1 ]}{\Rightarrow} \ 641 \ | \ F_5 $. \hfill $ \square $

\medskip

\underline{Proof VI} (Broda, see Chapter XV of \cite{Dickson})

$ \left[ 4 \cdot 5 \cdot 2^5 + 1 \ \textrm{is a prime number} \ \wedge \ 4 \cdot 5 \cdot 2^5 + 1 \ \nmid \ 2 > 0 \right] \ \overset{\tiny \begin{array}{c} \textrm{Fermat's} \\ \textrm{little} \\ \textrm{theorem} \end{array}}{\Rightarrow} \ 4 \cdot 5 \cdot 2^5 + 1 \ | \  ( 2^{2^5} )^{4 \cdot 5} - 1 $. \hfill $ [ 1 ] $

$ \left[ 4 \cdot 5 \cdot 2^5 + 1 \ \nmid \ ( 2^{2^5} - 1 ) \cdot ( ( 2^{2^5} )^2 + 1 ) \cdot \displaystyle\sum_{k = 0}^4 \left( ( 2^{2^5} )^{4 \cdot k} \right) \ \wedge \ ( 2^{2^5} )^4 \neq 1 \right] \ \overset{\tiny \begin{array}{c} \textrm{Sum of a} \\ \textrm{geometric} \\ \textrm{series} \end{array}}{\Rightarrow} $

$ 4 \cdot 5 \cdot 2^5 + 1 \ \nmid \ ( 2^{2^5} - 1 ) \cdot ( ( 2^{2^5} )^2 + 1 ) \cdot \dfrac{( 2^{2^5} )^{4 \cdot 5} - 1}{( 2^{2^5} )^4 - 1} \ \Rightarrow \ 4 \cdot 5 \cdot 2^5 + 1 \ \nmid \ \dfrac{( 2^{2^5} )^{4 \cdot 5} - 1}{2^{2^5} + 1} \ \overset{[ 1 ]}{\Rightarrow} \ \underbrace{4 \cdot 5 \cdot 2^5 + 1}_{641} \ | \ \underbrace{2^{2^5} + 1}_{F_5} $. \hfill $ \square $

\subsubsection{Generalization of the first proof}

\begin{prop}{\label{ZeroGen}}
    $ A \ | \ B \cdot C + E $, for every $ A , B , C , D , E \in \mathbb{Z} $ such that $ 0 \neq A \ | \ D + B , D \cdot C - E $.
\end{prop}

\underline{Proof}

Already seen in the example of Subsection \ref{SectionExplain}. \hfill $ \square $

\medskip

A possible particularization is the following one.

\begin{prop}{\label{FirstGen}}
    $ k \cdot l + m \ | \ q \cdot l^{2 \cdot r} + m^{2 \cdot r} $, for every $ k , l , m , q , r \in \mathbb{Z} $ such that $ r > 0 $ and $ k \cdot l + m \ | \ k^{2 \cdot r} + q $.
\end{prop}

\underline{Proof}

Let $ A = k \cdot l + m $, $ B = q $, $ C = l^{2 \cdot r} $, $ D = k^{2 \cdot r} $ and $ E = m^{2 \cdot r} $.

Then $ 0 \neq A \ | \ D + B , D \cdot C - E $.

Applying Proposition \ref{ZeroGen}, $ A \ | \ B \cdot C + E $; i.e., $ k \cdot l + m \ | \ q \cdot l^{2 \cdot r} + m^{2 \cdot r} $. \hfill $ \square $

\medskip

Particularizing again, the following result is obtained.

\begin{thm}{\label{Particularizacion}}
    $ k \cdot 2^s + 1 \ | \ F_n $, for every $ k , n , r , s \in \mathbb{N}^+ $ such that $ r \cdot s \leq 2^{n - 1} $ and $ k \cdot 2^s + 1 \ | \ k^{2 \cdot r} + 2^{2^n - 2 \cdot r \cdot s} $.
\end{thm}

\newpage

\underline{Proof}

Let $ l = 2^s $, $ m = 1 $ and $ q = 2^{2^n - 2 \cdot r \cdot s} $.

Then $ k \cdot l + m \ | \ k^{2 \cdot r} + q $.

Applying Proposition \ref{FirstGen}, $ k \cdot l + m \ | \ q \cdot l^{2 \cdot r} + m^{2 \cdot r} $; i.e., $ k \cdot 2^s + 1 \ | \ F_n $. \hfill $ \square $

\begin{problem}
    Do there exist $ k , n , r \in \mathbb{N}^+ $ such that $ k \cdot 2^{n + 2} + 1 $ is a prime factor of $ F_n $, $ r \cdot ( n + 2 ) \leq 2^{n - 1} $ and $ k \cdot 2^{n + 2} + 1 \ \nmid \ k^{2 \cdot r} + 2^{2^n - 2 \cdot r \cdot ( n + 2 )} $?
\end{problem}

The following Maple program computes the set $ \{ k \in \{ s , ... , t \} \ | \ k \cdot 2^{n + 2} + 1 \ | \ k^{2 \cdot r} + 2^{2^n - 2 \cdot r \cdot ( n + 2 )} \} $, with $ r $ being the natural number that is nearest to $ \dfrac{2^{n - 1}}{n + 3} $ $ \left( \textrm{note that} \ 2 \cdot \dfrac{2^{n - 1}}{n + 3} = 2^n - 2 \cdot \dfrac{2^{n - 1}}{n + 3} \cdot ( n + 2 ) \right) $, for given $ n , s , t \in \mathbb{N}^+ $ such that $ s \leq t $.

\begin{code}{}
test:=proc(n,s,t)
	local b, k, l, m, q, r, F:
	F:={}:
	b:=2^n:
	l:=4*b:
	m:=1:
	r:=round(b/(2*n+6)):
	q:=2^(b-2*r*(n+2)):
	for k from s to t do
		if (k^(2*r)+q) mod (k*l+m) = 0 then
			F:=`union`(F,{k}):
		fi:
	od:
	return F:
end proc:
\end{code}

For example, the following command returns $ \{ 7 , 1588 , 3892 \} $.

\begin{code}{}
test(12,1,5000);
\end{code}

\subsubsection{Generalization of the second proof}

\begin{prop}{\label{SecondGen}}
    $ b - c \ | \ \dfrac{g \cdot b^d}{b - c - g} + c^d $, for every $ b , c , d , g \in \mathbb{Z} $ such that $ d > 0 \neq b - c \neq g $.
\end{prop}

\underline{Proof}

$ b - c \equiv 0 \ ( \textrm{mod}. \ b - c ) \ \Rightarrow \ b - c - g \equiv - g \ ( \textrm{mod}. \ b - c ) \ \Rightarrow \ b^d \equiv - \dfrac{g \cdot b^d}{b - c - g} \ ( \textrm{mod}. \ b - c ) $. \hfill $ [ 1 ] $

$ b \equiv c \ ( \textrm{mod}. \ b - c ) \ \Rightarrow \ b^d \equiv c^d \ ( \textrm{mod}. \ b - c ) \ \overset{[ 1 ]}{\Rightarrow} \ - \dfrac{g \cdot b^d}{b - c - g} \equiv c^d \ ( \textrm{mod}. \ b - c ) \ \Rightarrow \ b - c \ | \ \dfrac{g \cdot b^d}{b - c - g} + c^d $. \hfill $ \square $

\medskip

This generalization and Proposition \ref{FirstGen} have the following common particularization.

\begin{prop}{\label{CommonPart}}
    $ k \cdot l + m \ | \ ( k \cdot l + m - k^{2 \cdot r} ) \cdot l^{2 \cdot r} + m^{2 \cdot r} $, for every $ k , l , m , r \in \mathbb{Z} $ such that $ r > 0 \neq k $.
\end{prop}

\underline{Proof I} (particularization of Proposition \ref{FirstGen})

Let $ q = k \cdot l + m - k^{2 \cdot r} $.

Then $ k \cdot l + m \ | \ k^{2 \cdot r} + q $.

Applying Proposition \ref{FirstGen}, $ k \cdot l + m \ | \ q \cdot l^{2 \cdot r} + m^{2 \cdot r} $; i.e., $ k \cdot l + m \ | \ ( k \cdot l + m - k^{2 \cdot r} ) \cdot l^{2 \cdot r} + m^{2 \cdot r} $. \hfill $ \square $

\medskip

\underline{Proof II} (particularization of Proposition \ref{SecondGen})

Let $ b = k \cdot l $, $ c = - m $, $ d = 2 \cdot r $ and $ g = k \cdot l + m - k^{2 \cdot r} $.

$ r > 0 \neq k $ implies that $ d > 0 \neq b - c \neq g $.

Applying Proposition \ref{SecondGen}, $ b - c \ | \ \dfrac{g \cdot b^d}{b - c - g} + c^d $; i.e., $ k \cdot l + m \ | \ ( k \cdot l + m - k^{2 \cdot r} ) \cdot l^{2 \cdot r} + m^{2 \cdot r} $. \hfill $ \square $

\medskip

Particularizing again, the following result is obtained.

\begin{thm}{\label{Particularizacion2}}
    $ s \cdot 2^t + 1 \ | \ F_n $, for every $ l , n , r , s , t \in \mathbb{N}^+ $ such that $ ( s \cdot 2^t + 1 ) \cdot l^{2 \cdot r} - ( s \cdot 2^t )^{2 \cdot r} + 1 = F_n $ and $ l \ | \ s \cdot 2^t $.
\end{thm}

\underline{Proof}

Let $ k = \dfrac{s}{l} \cdot 2^t $ and $ m = 1 $.

Then $ r > 0 \neq k $.

Applying Proposition \ref{CommonPart}, $ k \cdot l + m \ | \ ( k \cdot l + m - k^{2 \cdot r} ) \cdot l^{2 \cdot r} + m^{2 \cdot r} $; i.e.,

$ s \cdot 2^t + 1 \ | \ \left( s \cdot 2^t + 1 - \left( \dfrac{s}{l} \cdot 2^t \right)^{2 \cdot r} \right) \cdot l^{2 \cdot r} + 1 $.

And, by the hypothesis that $ ( s \cdot 2^t + 1 ) \cdot l^{2 \cdot r} - ( s \cdot 2^t )^{2 \cdot r} + 1 = F_n $, $ s \cdot 2^t + 1 \ | \ F_n $. \hfill $ \square $

\begin{problem}
    Determine the numbers $ n \in \mathbb{N}^+ $ for which there exist $ r , s , t \in \mathbb{N}^+ $ and a factor $ l $ of $ s \cdot 2^t $ such that $ ( s \cdot 2^t + 1 ) \cdot l^{2 \cdot r} - ( s \cdot 2^t )^{2 \cdot r} + 1 = F_n $.
\end{problem}

\subsubsection{Generalization of the third proof}

\begin{prop}[Baaz]{\label{PropBaaz}}
    $ v \cdot 2^u + 1 \ | \ 2^{2 \cdot u \cdot x} \cdot ( v \cdot 2^u + 1 - v^{2 \cdot x} ) + 1 $, for every $ u , v , x \in \mathbb{N}^+ $.
\end{prop}

\underline{Proof I} (direct generalization)

$ 2 \cdot x \ \textrm{is even} \ \Rightarrow \ \displaystyle\sum_{k = 0}^{2 \cdot x} \left( \displaystyle\binom{2 \cdot x}{k} \cdot ( v \cdot 2^u + 1 )^{2 \cdot x - k} \cdot ( - 1 )^k \right) = \displaystyle\sum_{k = 0}^{2 \cdot x - 1} \left( \displaystyle\binom{2 \cdot x}{k} \cdot ( v \cdot 2^u + 1 )^{2 \cdot x - k} \cdot ( - 1 )^k \right) + 1 $. \hfill $ [ 1 ] $

$ v^{2 \cdot x} = ( v \cdot 2^u + 1 ) - ( v \cdot 2^u + 1 - v^{2 \cdot x} ) \ \Rightarrow \ v^{2 \cdot x} \cdot 2^{u \cdot 2 \cdot x} = ( v \cdot 2^u + 1 ) \cdot 2^{u \cdot 2 \cdot x} - ( v \cdot 2^u + 1 - v^{2 \cdot x} ) \cdot 2^{u \cdot 2 \cdot x} $. \hfill $ [ 2 ] $
    
$ v \cdot 2^u = v \cdot 2^u + 1 + ( - 1 ) \ \Rightarrow \ v^{2 \cdot x} \cdot 2^{u \cdot 2 \cdot x} = ( v \cdot 2^u + 1 + ( - 1 ) )^{2 \cdot x} \ \overset{\tiny \begin{array}{c} \textrm{Binomial} \\ \textrm{theorem} \end{array}}{\Rightarrow} $
    
$ v^{2 \cdot x} \cdot 2^{u \cdot 2 \cdot x} = \displaystyle\sum_{k = 0}^{2 \cdot x} \left( \displaystyle\binom{2 \cdot x}{k} \cdot ( v \cdot 2^u + 1 )^{2 \cdot x - k} \cdot ( - 1 )^k \right) \ \overset{[ 1 ]}{\Rightarrow} $

$ v^{2 \cdot x} \cdot 2^{u \cdot 2 \cdot x} = \displaystyle\sum_{k = 0}^{2 \cdot x - 1} \left( \displaystyle\binom{2 \cdot x}{k} \cdot ( v \cdot 2^u + 1 )^{2 \cdot x - k} \cdot ( - 1 )^k \right) + 1 \ \Rightarrow $
    
$ v^{2 \cdot x} \cdot 2^{u \cdot 2 \cdot x} = ( v \cdot 2^u + 1 ) \cdot \displaystyle\sum_{k = 0}^{2 \cdot x - 1} \left( \displaystyle\binom{2 \cdot x}{k} \cdot ( v \cdot 2^u + 1 )^{2 \cdot x - k - 1} \cdot ( - 1 )^k \right) + 1 \ \overset{[ 2 ]}{\Rightarrow} $
    
$ ( v \cdot 2^u + 1 ) \cdot 2^{u \cdot 2 \cdot x} - ( v \cdot 2^u + 1 - v^{2 \cdot x} ) \cdot 2^{u \cdot 2 \cdot x} = ( v \cdot 2^u + 1 ) \cdot \displaystyle\sum_{k = 0}^{2 \cdot x - 1} \left( \displaystyle\binom{2 \cdot x}{k} \cdot ( v \cdot 2^u + 1 )^{2 \cdot x - k - 1} \cdot ( - 1 )^k \right) + 1 \ \Rightarrow $
    
$ ( v \cdot 2^u + 1 ) \cdot \left( 2^{u \cdot 2 \cdot x} - \displaystyle\sum_{k = 0}^{2 \cdot x - 1} \left( \displaystyle\binom{2 \cdot x}{k} \cdot ( v \cdot 2^u + 1 )^{2 \cdot x - k - 1} \cdot ( - 1 )^k \right) \right) = ( v \cdot 2^u + 1 - v^{2 \cdot x} ) \cdot 2^{u \cdot 2 \cdot x} + 1 \ \Rightarrow $

$ v \cdot 2^u + 1 \ | \ ( v \cdot 2^u + 1 - v^{2 \cdot x} ) \cdot 2^{u \cdot 2 \cdot x} + 1 $. \hfill $ \square $

\medskip

\underline{Proof II} (particularization of Proposition \ref{CommonPart})

Let $ k = v $, $ l = 2^u $, $ m = 1 $ and $ r = x $.

Then $ r > 0 \neq k $.

Applying Proposition \ref{CommonPart}, $ k \cdot l + m \ | \ ( k \cdot l + m - k^{2 \cdot r} ) \cdot l^{2 \cdot r} + m^{2 \cdot r} $; i.e., $ v \cdot 2^u + 1 \ | \ 2^{2 \cdot u \cdot x} \cdot ( v \cdot 2^u + 1 - v^{2 \cdot x} ) + 1 $. \hfill $ \square $

\medskip

The next particularization follows immediately.

\begin{thm}[Baaz]
    $ v \cdot 2^u + 1 \ | \ F_n $, for every $ n , u , v , x \in \mathbb{N}^+ $ such that $ 2^{2 \cdot u \cdot x} \cdot ( v \cdot 2^u + 1 - v^{2 \cdot x} ) + 1 = F_n $.
\end{thm}

\underline{Proof I} (particularization of Proposition \ref{PropBaaz})

Immediate. \hfill $ \square $

\underline{Proof II} (particularization of Theorem \ref{Particularizacion2})

Let $ l = 2^u $, $ r = x $, $ s = v $ and $ t = u $.

Then $ ( s \cdot 2^t + 1 ) \cdot l^{2 \cdot r} - ( s \cdot 2^t )^{2 \cdot r} + 1 = F_n $ and $ l \ | \ s \cdot 2^t $.

Applying Theorem \ref{Particularizacion2}, $ s \cdot 2^t + 1 \ | \ F_n $; i.e., $ v \cdot 2^u + 1 \ | \ F_n $. \hfill $ \square $

\subsubsection{Generalization of the fourth proof}

\begin{prop}{\label{ThirdGen}}
    $ a^h + i^4 \ | \ a^{h + 4 \cdot f} + ( a^h - i \cdot a^f + i^4 )^4 $, for every $ a , f , h , i \in \mathbb{Z} $ such that $ a^h + i^4 \neq 0 < f , h $.
\end{prop}

\underline{Proof I} (direct generalization)

$ a^f - i^3 - ( a^f - i^3 ) = 0 \ \Rightarrow \ - i \cdot ( a^f - i^3 ) + i \cdot ( a^f - i^3 ) = 0 \ \Rightarrow \ ( a^h - i \cdot a^f + i^4 ) + a^f \cdot i - i^4 = a^h $. \hfill $ [ 1 ] $

$ h + 4 \cdot f = h + 4 \cdot f \ \Rightarrow \ a^{h + 4 \cdot f} = a^{h + 4 \cdot f} \ \Rightarrow \ a^{h + 4 \cdot f} + ( a^h - i \cdot a^f + i^4 )^4 = a^{h + 4 \cdot f} + ( a^h - i \cdot a^f + i^4 )^4 \ \overset{[ 1 ]}{\Rightarrow} $

$ a^{h + 4 \cdot f} + ( a^h - i \cdot a^f + i^4 )^4 = ( ( a^h - i \cdot a^f + i^4 ) + a^f \cdot i - i^4 ) \cdot a^{4 \cdot f} + ( a^h - i \cdot a^f + i^4 )^4 \ \Rightarrow $

$ a^{h + 4 \cdot f} + ( a^h - i \cdot a^f + i^4 )^4 = ( a^h + i^4 ) \cdot a^{4 \cdot f} + ( ( a^h - i \cdot a^f + i^4 )^4 - i^4 \cdot a^{4 \cdot f} ) \ \Rightarrow $

$ a^{h + 4 \cdot f} + ( a^h - i \cdot a^f + i^4 )^4 = $

$ ( a^h + i^4 ) \cdot a^{4 \cdot f} + ( ( a^h - i \cdot a^f + i^4 )^2 + i^2 \cdot a^{2 \cdot f} ) \cdot ( ( a^h - i \cdot a^f + i^4 ) + i \cdot a^f ) \cdot ( ( a^h - i \cdot a^f + i^4 ) - i \cdot a^f ) \ \Rightarrow $

$ a^{h + 4 \cdot f} + ( a^h - i \cdot a^f + i^4 )^4 = ( a^h + i^4 ) \cdot ( a^{4 \cdot f} + ( ( a^h - i \cdot a^f + i^4 )^2 + i^2 \cdot a^{2 \cdot f} ) \cdot ( ( a^h - i \cdot a^f + i^4 ) - i \cdot a^f ) ) \ \Rightarrow $

$ a^h + i^4 \ | \ a^{h + 4 \cdot f} + ( a^h - i \cdot a^f + i^4 )^4 $. \hfill $ \square $

\medskip

\underline{Proof II} (particularization of Proposition \ref{FirstGen})

Let $ k = i $, $ l = a^f $, $ m = a^h - i \cdot a^f + i^4 $, $ q = a^h $, $ r = 2 $.

Then $ r > 0 $ and $ k \cdot l + m \ | \ k^{2 \cdot r} + q $.

Applying Proposition \ref{FirstGen}, $ k \cdot l + m \ | \ q \cdot l^{2 \cdot r} + m^{2 \cdot r} $; i.e., $ a^h + i^4 \ | \ a^{h + 4 \cdot f} + ( a^h - i \cdot a^f + i^4 )^4 $. \hfill $ \square $

\medskip

This generalization and Proposition \ref{SecondGen} have the following common particularization.

\begin{prop}{\label{CommonPart2}}
    $ a^h + i^4 \ | \ a^{h + 4 \cdot f} + ( a^h - i \cdot a^f + i^4 )^4 $, for every $ a , f , h , i \in \mathbb{Z} $ such that $ a^h + i^4 \neq 0 < f , h $ and $ i \neq 0 $.
\end{prop}

\underline{Proof I} (particularization of Proposition \ref{ThirdGen})

Immediate. \hfill $ \square $

\medskip

\underline{Proof II} (particularization of Proposition \ref{SecondGen})

Let $ b = i \cdot a^f $, $ c = i \cdot a^f - a^h - i^4 $, $ d = 4 $ and $ g = a^h $.

$ 0 \neq i^4 \neq - a^h $ implies that $ d > 0 \neq b - c = a^h + i^4 \neq g = a^h $.

Applying Proposition \ref{SecondGen}, $ b - c \ | \ \dfrac{g \cdot b^d}{b - c - g} + c^d $; i.e., $ a^h + i^4 \ | \ a^{h + 4 \cdot f} + ( a^h - i \cdot a^f + i^4 )^4 $. \hfill $ \square $

The following particularization is easily derived.

\begin{thm}{\label{Particularizacion3}}
    $ 2^{2^n - 4 \cdot ( n + 2 )} + i^4 \ | \ F_n $, for every $ i , n \in \mathbb{N}^+ $ such that $ n > 4 $ and $ i \cdot 2^{n + 2} + 1 = 2^{2^n - 4 \cdot ( n + 2 )} + i^4 $.
\end{thm}

\underline{Proof}

Let $ a = 2 $, $ f = n + 2 $ and $ h = 2^n - 4 \cdot ( n + 2 ) $.

$ n > 4 $ implies that $ a^h + i^4 \neq 0 < f , h $ and $ i \neq 0 $.

Applying Proposition \ref{CommonPart2}, $ a^h + i^4 \ | \ a^{h + 4 \cdot f} + ( a^h - i \cdot a^f + i^4 )^4 $; i.e.

$ 2^{2^n - 4 \cdot ( n + 2 )} + i^4 \ | \ 2^{2^n} + ( 2^{2^n - 4 \cdot ( n + 2 )} - i \cdot 2^{n + 2} + i^4 )^4 $.

And, by the hypothesis that $ i \cdot 2^{n + 2} + 1 = 2^{2^n - 4 \cdot ( n + 2 )} + i^4 $, $ 2^{2^n - 4 \cdot ( n + 2 )} + i^4 \ | \ F_n $. \hfill $ \square $

\begin{problem}
    Determine the numbers $ n \in \mathbb{N}^+ $ for which there exists $ i \in \mathbb{N}^+ $ such that
    
    $ i \cdot 2^{n + 2} + 1 = 2^{2^n - 4 \cdot ( n + 2 )} + i^4 $.
\end{problem}

\subsubsection{Generalization of the fifth proof}

\begin{prop}{\label{FourthGen}}
    $ i \ | \ ( a + k \cdot c )^d - ( a + b \cdot c )^d + g \cdot h $, for every $ a , b , c , d , g , h , i , k \in \mathbb{Z} $ such that $ d > 0 \neq i \ | \ k - b , g $.
\end{prop}

\underline{Proof}

$ i \ | \ g \ \Rightarrow \ i \ | \ g \cdot h \ \Rightarrow \ i \ | \ ( a + b \cdot c )^d - ( ( a + b \cdot c )^d - g \cdot h ) \ \Rightarrow $

$ i \ | \ ( ( a + b \cdot c )^d - ( a + k \cdot c )^d ) + ( - ( ( a + b \cdot c )^d - g \cdot h ) + ( a + k \cdot c )^d ) $. \hfill $ [ 1 ] $

$ i \ | \ k - b \ \Rightarrow \ i \ | \ k \cdot c - b \cdot c \ \Rightarrow \ i \ | \ ( a + k \cdot c ) - ( a + b \cdot c ) \ \Rightarrow $

$ i \ | \ ( a + k \cdot c )^d - ( a + b \cdot c )^d \ \Rightarrow \ i \ | \ ( a + b \cdot c )^d - ( a + k \cdot c )^d \ \overset{[ 1 ]}{\Rightarrow} \ i \ | \ ( a + k \cdot c )^d - ( a + b \cdot c )^d + g \cdot h $. \hfill $ \square $

\medskip

The following proposition is a possible particularization.

\begin{thm}{\label{Particularizacion4}}
    $ i \ | \ F_n $, for every $ c , i , n \in \mathbb{N}^+ $ such that $ i \ | \ ( 2^{2^{n - 1}} - i \cdot c )^2 + 1 $.
\end{thm}

\underline{Proof}

Let $ h \in \mathbb{N}^+ $ such that $ ( 2^{2^{n - 1}} - i \cdot c )^2 = i \cdot h - 1 $, $ a = 2^{2^{n - 1}} - ( i - 1 ) \cdot c $, $ b = - 1 $, $ d = 2 $, $ g = i $ and $ k = i - 1 $.

Then $ d > 0 \neq i \ | \ k - b , g $.

Applying Proposition \ref{FourthGen}, $ i \ | \ ( a + k \cdot c )^d - ( a + b \cdot c )^d + g \cdot h = 2^{2^n} - ( 2^{2^{n - 1}} - i \cdot c )^2 + i \cdot h = F_n $. \hfill $ \square $

\begin{problem}
    Determine the numbers $ n \in \mathbb{N}^+ $ for which there exist $ c , i \in \mathbb{N}^+ $ such that $ i \ | \ ( 2^{2^{n - 1}} - i \cdot c )^2 + 1 $.
\end{problem}

\subsubsection{Generalization of the sixth proof}{\label{Little}}

\begin{prop}{\label{Broda}}
    $ ( C + 1 ) \cdot A + 1 \ | \ D $, for every $ A , B , C , D \in \mathbb{N} $ such that $ B > 1 $, $ D \ | \ B^A - 1 $ and $ ( C + 1 ) \cdot A + 1 $ is a prime number that does not divide $ B \cdot \dfrac{B^A - 1}{D} \cdot \displaystyle\sum_{k = 0}^C \left( ( B^A )^k \right) $.
\end{prop}

\underline{Proof}

Note that the hypothesis of $ ( C + 1 ) \cdot A + 1 $ being a prime number yields that $ A > 0 $; and the hypothesis that it does not divide $ B \cdot \dfrac{B^A - 1}{D} \cdot \displaystyle\sum_{k = 0}^C \left( ( B^A )^k \right) $ implies that it does not divide $ B $ and neither $ \dfrac{B^A - 1}{D} \cdot \displaystyle\sum_{k = 0}^C \left( ( B^A )^k \right) $.

$ \left[ ( C + 1 ) \cdot A + 1 \ \textrm{is a prime number} \ \wedge \ ( C + 1 ) \cdot A + 1 \ \nmid \ B > 0 \right] \ \overset{\tiny \begin{array}{c} \textrm{Fermat's} \\ \textrm{little} \\ \textrm{theorem} \end{array}}{\Rightarrow} \ ( C + 1 ) \cdot A + 1 \ | \  B^{( C + 1 ) \cdot A} - 1 $. \hfill $ [ 1 ] $

$ \left[ ( C + 1 ) \cdot A + 1 \ \nmid \ \dfrac{B^A - 1}{D} \cdot \displaystyle\sum_{k = 0}^C \left( ( B^A )^k \right) \ \wedge \ B^A \neq 1 \right] \ \overset{\tiny \begin{array}{c} \textrm{Sum of a} \\ \textrm{geometric} \\ \textrm{series} \end{array}}{\Rightarrow} $

$ ( C + 1 ) \cdot A + 1 \ \nmid \ \dfrac{B^A - 1}{D} \cdot \dfrac{( B^A )^{C + 1} - 1}{B^A - 1} \ \Rightarrow \ ( C + 1 ) \cdot A + 1 \ \nmid \ \dfrac{( B^A )^{C + 1} - 1}{D} \ \overset{[ 1 ]}{\Rightarrow} \ ( C + 1 ) \cdot A + 1 \ | \ D $. \hfill $ \square $

\medskip

Two possible particularizations are the following ones.

\begin{thm}{\label{Particularizacion5}}
    $ m \cdot 2^{n + 2} + 1 \ | \ F_n $, for every $ m , n \in \mathbb{N}^+ $ such that $ m \cdot 2^{n + 2} + 1 $ is a prime number that does not divide $ \dfrac{( F_{n + 2} - 1 )^m - 1}{F_n} $.
\end{thm}

\underline{Proof}

Let $ A = 2^{n + 2} $, $ B = 2 $, $ C = m - 1 $ and $ D = F_n $.

Then $ B > 1 $ and $ ( C + 1 ) \cdot A + 1 $ is a prime number that does not divide $ 2 \cdot \dfrac{( F_{n + 2} - 1 )^m - 1}{F_n} $;

i.e., $ B \cdot \dfrac{B^A - 1}{D} \cdot \displaystyle\sum_{k = 0}^C \left( ( B^A )^k \right) $.

In addition, applying Theorem \ref{RecDef}, $ D \ | \ B^A - 1 = \displaystyle\prod_{k = 0}^{n + 1} ( F_k ) $.

Applying Proposition \ref{Broda}, $ ( C + 1 ) \cdot A + 1 \ | \ D $; i.e., $ m \cdot 2^{n + 2} + 1 \ | \ F_n $. \hfill $ \square $

\begin{problem}
    Do there exist $ m , n \in \mathbb{N}^+ $ such that $ m \cdot 2^{n + 2} + 1 $ is a prime number that divides both $ F_n $ and $ \dfrac{( F_{n + 2} - 1)^m - 1}{F_n} $?
\end{problem}

\begin{thm}
    $ m \cdot 2^{n + 2} + 1 \ | \ \displaystyle\prod_{k = 0}^{n + 1} ( F_k ) $, for every $ m , n \in \mathbb{N}^+ $ such that $ m \cdot 2^{n + 2} + 1 $ is a prime number that does not divide $ \dfrac{( F_{n + 2} - 1 )^m - 1}{F_{n + 2} - 2} $.
\end{thm}

\underline{Proof}

Let $ A = 2^{n + 2} $, $ B = 2 $, $ C = m - 1 $ and $ D = F_{n + 2} - 2 $.

Then $ B > 1 $, $ D \ | \ B^A - 1 $ and $ ( C + 1 ) \cdot A + 1 $ is a prime number that does not divide $ 2 \cdot \dfrac{( F_{n + 2} - 1 )^m - 1}{F_{n + 2} - 2} $; i.e., $ B \cdot \dfrac{B^A - 1}{D} \cdot \displaystyle\sum_{k = 0}^C \left( ( B^A )^k \right) $.

Applying Proposition \ref{Broda}, $ ( C + 1 ) \cdot A + 1 \ | \ D $; i.e., $ m \cdot 2^{n + 2} + 1 \ | \ F_{n + 2} - 2 $.

And, by Theorem \ref{RecDef}, $ m \cdot 2^{n + 2} + 1 \ | \ \displaystyle\prod_{k = 0}^{n + 1} ( F_k ) $. \hfill $ \square $

\begin{problem}
    Do there exist $ m , n \in \mathbb{N}^+ $ such that $ m \cdot 2^{n + 2} + 1 $ is a prime number that divides both $ \displaystyle\prod_{k = 0}^{n + 1} ( F_k ) $ and $ \dfrac{( F_{n + 2} - 1)^m - 1}{F_{n + 2} - 2} $?
\end{problem}

For example, $ 37 \cdot 2^{14 + 2} + 1 $, which is a prime factor of $ F_9 $, does not divide $ \dfrac{( F_{14 + 2} - 1)^{37} - 1}{F_{14 + 2} - 2} $ but divides $ \dfrac{( F_{14 + 2} - 1)^{37} - 1}{F_{14}} $.

During the rest of this subsubsection, the function from $ \mathbb{N}^+ \times \mathbb{N}^+ $ to $ \{ 0 , 1 \} $ that, for every $ m , n \in \mathbb{N}^+ $, maps $ ( m , n ) $ into $ 1 $ if, and only if, $ m \cdot 2^{n + 2} + 1 $ is a prime number that does not divide $ \dfrac{( F_{n + 2} - 1)^m - 1}{F_{n + 2} - 2} $ $ \left( \textrm{resp}. \ \dfrac{( F_{n + 2} - 1)^m - 1}{F_{n}} \right) $, will be denoted by $ A $ (resp. $ B $).

The following Maple function compute the images of $ A $ and $ B $, returning ``true'' if the image is $ 1 $ and ``false'' otherwise.

\begin{code}{}
A:=(m,n)->isprime(m*2^(n+2)+1) and not ((2^(m*2^(n+2))-1)/(2^(2^(n+2))-1)) mod (m*2^(n+2)+1) = 0:
B:=(m,n)->isprime(m*2^(n+2)+1) and not ((2^(m*2^(n+2))-1)/(2^(2^n)+1)) mod (m*2^(n+2)+1) = 0:
\end{code}    

For example, the following command returns ``true''.

\begin{code}{}
A(5,5);
\end{code}

\begin{problem}
    If $ n \in \mathbb{N}^+ $, what is the cardinality of $ \{ m \in \mathbb{N}^+ \ | \ A ( m , n ) = 1 \} $?
\end{problem}

For example, $ A ( 8 , 11 ) = A ( 14 , 11 ) = A ( 39 , 11 ) = A ( 119 , 11 ) = 1 $.

\begin{problem}
    Assuming that $ \{ m \in \mathbb{N}^+ \ | \ A ( m , n ) = 1 \} $ is nonempty, for every $ n \in \mathbb{N}^+ $, find, if possible, an explicit formula for the sequence whose $ n^\textrm{th} $ term is $ \textrm{min} ( \{ m \in \mathbb{N}^+ \ | \ A ( m , n + 1 ) = 1 \} ) $, for every $ n \in \mathbb{N} $.
\end{problem}

The first 14 terms of the previous problem's sequence are $ 2 , 1 , 8 , 4 , 2 , 1 , 128 , 64 , 32 , 16 , 8 , 4 , 2 , 1 $.

\begin{problem}
    If $ m \in \mathbb{N}^+ $, what is the cardinality of $ \{ n \in \mathbb{N}^+ \ | \ A ( m , n ) = 1 \} $?
\end{problem}

For example, $ A ( 2 , 1 ) = A ( 2 , 5 ) = A ( 2 , 13 ) = 1 $.

\begin{problem}
    If $ m , n \in \mathbb{N}^+ $ are such that $ m $ is odd and $ A ( m , n ) = 1 $, what is the cardinality of
    
    $ \{ k \in \mathbb{N} \ | \ A ( 2^i \cdot m , n - i ) , \ \textrm{for every} \ i \in \{ 0 , ... , k \} \} $?
\end{problem}

For example, $ A ( 2^0 \cdot 37 , 14 - 0 ) = ... = A ( 2^6 \cdot 37 , 14 - 6 ) = 1 $.

\begin{problem}
    What is the cardinality of $ \{ n \in \mathbb{N}^+ \ | \ A ( n , n ) = 1 \} $?
\end{problem}

For example, $ A ( 4 , 4 ) = A ( 5 , 5 ) = 1 $.

\begin{problem}
   What is the cardinality of $ \mathbb{N}^+ \backslash \{ n \in \mathbb{N}^+ \ | \ A ( 1 , n ) = A ( 2 , n - 1 ) = ... = A ( n , 1 ) = 0 \} $?
\end{problem}

For example, $ A ( 1 , 3 ) = A ( 2 , 2 ) = A ( 3 , 1 ) = 0 $.

The following Maple program (which requires to have loaded the function $ A $) computes the elements of the previous problem's set that are upper-bounded by a given $ r \in \mathbb{N} $.

\begin{code}{}
test:=proc(r)
	local k, n, F:
	F:={}:
	for n from 1 to r do
		for k from 0 to n do
			if A(k+1,n-k) then
				F:=`union`(F,{n}):
				break
			fi:
		od:
	od:
	return F:
end proc:
\end{code}

For example, the following command returns $ \{ 2 , 6 , 7 , 9 , 10 , 13 , 14 , 15 , 17 , 18 , 24 , 25 , 27 , 30 \} $.

\begin{code}{}
test(30);
\end{code}

\begin{problem}
   Find, if possible, an explicit formula for the strictly increasing integer sequence whose image is the previous problem's set, if it is infinite.
\end{problem}

\begin{problem}
   What is the cardinality of $ \mathbb{N}^+ \backslash \{ n \in \mathbb{N}^+ \ | \ B ( 1 , n ) = B ( 2 , n - 1 ) = ... = B ( n , 1 ) = 0 \} $?
\end{problem}

\begin{problem}
   Find, if possible, an explicit formula for the strictly increasing integer sequence whose image is the previous problem's set, if it is infinite.
\end{problem}

The first 6 terms of the previous problem's integer sequence are $ 2, 9, 10, 18, 27 , 30 $ (they can be calculated with the previous Maple program, by replacing 'A' by 'B').

\begin{problem}
    What is the cardinality of $ \{ ( m , n ) \in ( \mathbb{N}^+ )^2 \ | \ A ( m , n ) = A ( m + 1 , n - 1 ) = 1 \} $?
\end{problem}

For example, $ A ( 1 , 2 ) = A ( 2 , 1 ) = A ( 1 , 6 ) = A ( 2 , 5 ) = A ( 1 , 14 ) = A ( 2 , 13 ) = A ( 7 , 12 ) = A ( 8 , 11 ) = 1 $.

\begin{problem}
    If $ S $ is the nonempty subset of $ \mathbb{N} $ such that, for every $ s \in \mathbb{N} $, $ s \in S $ if, and only if, there exist $ m , n \in \mathbb{N}^+ $ such that $ A ( m + 0 , n - 0 ) = ... = A ( m + s , n - s ) = 1 $, is $ \textrm{max} ( S ) = 1 $?
\end{problem}

\begin{problem}
    If $ a \in \mathbb{N}^+ $, do there exist $ m , n \in \mathbb{N}^+ $ such that, for every $ i \in \{ m , ... , m + a \} $ and
    
    $ j \in \{ n , ... , n + a \} $, $ A ( i , j ) = 0 $?
\end{problem}

More related problems are described in Subsection \ref{SubsectionInteresting}.

\section{Denumerable simple graphs as integer sequences}

In this section a technique for translating denumerable simple graphs (i.e., denumerable graphs without loops) into certain integer sequences (and viceversa) is introduced. This translation can help to obtain arithmetical properties of denumerable simple graphs, and graphical properties of certain integer sequences.

\subsection{The bijection}

Roughly speaking, the idea consists of expressing denumerable simple graphs as ``limits'' of sequences of finite simple graphs that add exactly one vertex for each term, and reading the last row (except the last entry) of each corresponding adjacency matrix as a binary number, defining therefore an integer sequence.

Formally, if $ G $ is a graph sequence, then $ G $ is \textbf{vertex-by-vertex increasing} if, and only if, for every $ n \in \mathbb{N} $, $ G_n $ is simple, $ V ( G_n ) = \{ 1 , ... , n + 1 \} $ and $ E ( G_n ) \subseteq E ( G_{n + 1} ) $ (recall that $ V ( G_n ) $ and $ E ( G_n ) $ respectively denote the vertex set and the edge set of $ G_n $).

Note that, if $ G $ is a vertex-by-vertex increasing graph sequence, then $ E ( G_0 ) = \varnothing $.

From now on, given a vertex-by-vertex increasing graph sequence $ G $, the graph whose vertex set is $\mathbb{N}^+ $ and whose edge set is $ \displaystyle\bigcup_{n \in \mathbb{N}} ( E ( G_n ) ) $ will be denoted by $ G_\infty $.

In addition, from now on, considering:
    \begin{itemize}
        \item $ S $ to be the set of vertex-by-vertex increasing graph sequences,
        \item $ T = \{ s \in \mathbb{N}^\mathbb{N} \ | \ s_n < 2^{n + 1} , \ \textrm{for every} \ n \in \mathbb{N} \} $,
        \item $ f : S \longrightarrow T $ such that $ ( f ( G ) )_n = \displaystyle\sum_{k = 1}^{n + 1} \left( \dfrac{( \textrm{adj} ( G_{n + 1} ) )_{n + 2 , k}}{2^{k - n - 1}} \right) $, for every $ G \in S $ and $ n \in \mathbb{N} $ (recall that $ \textrm{adj} ( G_{n+1} ) $ denotes the adjacency matrix of $ G_{n+1} $),
    \end{itemize}

$ f $ will be denoted by $ \Phi $.

\begin{thm}
    $ \Phi $ is bijective.
\end{thm}

\newpage

\underline{Proof}

Given a vertex-by-vertex increasing graph sequence $ G $, let $ A $ be the sequence of hollow symmetric binary matrices (i.e., symmetric binary matrices whose main diagonal has only zeros) such that $ A_n = \textrm{adj} ( G_n ) $ (which is of order $ n + 1 $), for every $ n \in \mathbb{N} $, and $ s \in \mathbb{N}^\mathbb{N} $ such that $ s_n = \displaystyle\sum_{k = 1}^{n + 1} \left( \dfrac{( A_{n + 1} )_{n + 2 , k}}{2^{k - n - 1}} \right) $, for every $ n \in \mathbb{N} $. Then $ s_n \leq \displaystyle\sum_{k = 1}^{n + 1} \left( \dfrac{1}{2^{k - n - 1}} \right) = 2^{n + 1} - 1 < 2^{n + 1} $, for every $ n \in \mathbb{N} $.

Conversely, given $ s \in \mathbb{N}^\mathbb{N} $ such that $ s_n < 2^{n + 1} $, for every $ n \in \mathbb{N} $, the length of the binary representation of $ s_n $ is shorter than $ n + 2 $, for every $ n \in \mathbb{N} $. Hence, there exists a unique sequence $ A $ of hollow symmetric binary matrices such that $ A_n $ is of order $ n + 1 $ and $ s_n = \displaystyle\sum_{k = 1}^{n + 1} \left( \dfrac{( A_{n + 1} )_{n + 2 , k}}{2^{k - n - 1}} \right) $, for every $ n \in \mathbb{N} $. Therefore, there exists a unique vertex-by-vertex increasing graph sequence $ G $ such that $ \textrm{adj} ( G_n ) = A_n $, for every $ n \in \mathbb{N} $. \hfill $ \square $

\medskip

The following Maple program computes $ ( \Phi ( G ) )_n $, for a given $ \textrm{adj} ( ( \Phi ( G ) )_n ) $ (of certain vertex-by-vertex increasing graph sequence $ G $ and $ n \in \mathbb{N} $).

\begin{code}{}
with(LinearAlgebra):
Phi:=proc(M)
	local n:
	n:=[Dimension(M)][1]-2:
	return add(M[n+2,k]/2^(k-n-1),k=1..n+1):
end proc:
\end{code}

For example, the following command returns $ 7 $ (recall that $ 111 $ is the binary representation of $ 7 $).

\begin{code}{}
Phi(<<0,1,1,1>|<1,0,1,1>|<1,1,0,1>|<1,1,1,0>>);
\end{code}

And the following Maple program computes $ \textrm{adj} ( ( \Phi^{- 1} ( s ) )_n ) $, for given $ n \in \mathbb{N} $ and $ s \in \mathbb{N}^\mathbb{N} $ such that $ s_k < 2^{k + 1} $, for every $ k \in \mathbb{N} $.

\begin{code}{}
InvPhi:=proc(n,s)
	local i, j, t, M:
	t:=unapply(s,k):
	M:=Matrix(1..n+2,1..n+2,shape=symmetric):
	for i from 0 to n do
		for j from 1 to length(convert(t(i),binary)) do:
			M[i+2,i+2-j]:=convert(t(i),base,2)[j]:
		od:
	od:
	return M:
end proc:
\end{code}    

For example, the following command returns $ \left[ \begin{array}{cccc} 0 & 1 & 1 & 1 \\ 1 & 0 & 1 & 1 \\ 1 & 1 & 0 & 1 \\ 1 & 1 & 1 & 0 \end{array} \right] $.

\begin{code}{}
InvPhi(2,2^(k+1)-1);
\end{code}    

\subsection{Some basic examples}

The first immediate example shows the integer sequence that is associated to the empty graph whose vertex set is $ \mathbb{N}^+ $.

\begin{prop}
    If $ G \in \textrm{dom} ( \Phi ) $ is such that $ E ( G_n ) = \varnothing $, for every $ n \in \mathbb{N} $, then $ ( \Phi ( G ) )_n = 0 $, for every $ n \in \mathbb{N} $.
\end{prop}

\begin{center}
    \begin{tabular}{|c|c|c|}
    \hline
    $ G_3 $ & $ \textrm{adj} ( G_3 ) $ & \\
    \hline
    \raisebox{-18pt}{\begin{tikzpicture}[]
  \node at (0,0) (A) {$ 1 $};
  \node at (0,1) (B) {$ 2 $};
  \node at (1,0) (C) {$ 3 $};
  \node at (1,1) (D) {$ 4 $};
  \path[] ;
    \end{tikzpicture}} & $ \left[ \begin{array}{cccc} * & * & * & * \\ 0 & * & * & * \\ 0 & 0 & * & * \\ 0 & 0 & 0 & * \end{array} \right] $ & $ \begin{array}{l} ( \Phi ( G ) )_0 = 0 \\ ( \Phi ( G ) )_1 = 0 \\ ( \Phi ( G ) )_2 = 0 \end{array} $ \\
    \hline
    \end{tabular}
\end{center}

The following example, showing the integer sequence that is associated to the complete graph whose vertex set is $ \mathbb{N}^+ $, is also immediate.

\begin{prop}
    If $ G \in \textrm{dom} ( \Phi ) $ is such that $ E ( G_n ) = \{ ( i , j ) \in V ( G_n )^2 \ | \ i \neq j \} $, for every $ n \in \mathbb{N} $, then $ ( \Phi ( G ) )_n = 2^{n + 1} - 1 $, for every $ n \in \mathbb{N} $.
\end{prop}

\begin{center}
    \begin{tabular}{|c|c|c|}
    \hline
    $ G_3 $ & $ \textrm{adj} ( G_3 ) $ & \\
    \hline
    \raisebox{-18pt}{\begin{tikzpicture}[]
  \node at (0,0) (A) {$ 1 $};
  \node at (0,1) (B) {$ 2 $};
  \node at (1,0) (C) {$ 3 $};
  \node at (1,1) (D) {$ 4 $};
  \path[]
        (A) edge node {} (B)
		edge node {} (C)
		edge node {} (D)
		(B) edge node {} (C)
        edge node {} (D)
        (C) edge node {} (D);
    \end{tikzpicture}} & $ \left[ \begin{array}{cccc} * & * & * & * \\ 1 & * & * & * \\ 1 & 1 & * & * \\ 1 & 1 & 1 & * \end{array} \right] $ & $ \begin{array}{l} ( \Phi ( G ) )_0 = 1 \\ ( \Phi ( G ) )_1 = 3 \\ ( \Phi ( G ) )_2 = 7 \end{array} $ \\
    \hline
    \end{tabular}
\end{center}

And again immediate is the next example, showing the integer sequence that is associated to the star graph whose vertex set is $ \mathbb{N}^+ $.

\begin{prop}
    If $ G \in \textrm{dom} ( \Phi ) $ is such that $ E ( G_n ) = \{ ( i , j ) \in V ( G_n )^2 \ | \ 1 \in \{ i , j \} \ \textrm{and} \ i \neq j \} $, for every $ n \in \mathbb{N} $, then $ ( \Phi ( G ) )_n = 2^n $, for every $ n \in \mathbb{N} $.
\end{prop}

\begin{center}
    \begin{tabular}{|c|c|c|}
    \hline
    $ G_3 $ & $ \textrm{adj} ( G_3 ) $ & \\
    \hline
    \raisebox{-18pt}{\begin{tikzpicture}[]
  \node at (1,1) (A) {$ 1 $};
  \node at (0,0) (B) {$ 2 $};
  \node at (1,0) (C) {$ 3 $};
  \node at (2,0) (D) {$ 4 $};
  \path[]
        (A) edge node {} (B)
		edge node {} (C)
		edge node {} (D);
    \end{tikzpicture}} & $ \left[ \begin{array}{cccc} * & * & * & * \\ 1 & * & * & * \\ 1 & 0 & * & * \\ 1 & 0 & 0 & * \end{array} \right] $ & $ \begin{array}{l} ( \Phi ( G ) )_0 = 1 \\ ( \Phi ( G ) )_1 = 2 \\ ( \Phi ( G ) )_2 = 4 \end{array} $ \\
    \hline
    \end{tabular}
\end{center}

The next example is not so easy though. Consider $ r \in \mathbb{N}^+ $. It shows the integer sequence that is associated to the tree whose vertex set is $ \mathbb{N}^+ $ and in which each vertex has exactly $ r $ successors. In order to explain it, the following functions from $ \mathbb{N}^+ $ to $ \mathbb{N}^+ $ are defined.
\begin{itemize}
    \item $ l_r $ (the ``level'' function), which maps each positive natural number $ n $ into itself if $ r = 1 $ and to $ \left\lfloor \log_r ( ( r - 1 ) \cdot ( n - 1 ) + 1 ) \right\rfloor + 1 $ otherwise. Some examples are the following ones.
    \begin{itemize}
        \item $ ( l_2 ( 1 ) , ... , l_2 ( 4 ) ) = ( 1 , 2 , 3 , 4 ) $,
        \item $ ( l_2 ( 1 ) , ... , l_2 ( 8 ) ) = ( 1 , 2 , 2 , 3 , 3 , 3 , 3 , 4 ) $,
        \item $ ( l_3 ( 1 ) , ... , l_3 ( 14 ) ) = ( 1 , 2 , 2 , 2 , 3 , 3 , 3 , 3 , 3 , 3 , 3 , 3 , 3 , 4 ) $.
    \end{itemize}
    \item $ m_r $, which maps each positive natural number $ n $ into $ \textrm{min} ( \{ i \in \mathbb{N}^+ \ | \ l_r ( i ) = l_r ( n ) \} ) $ (minimum number that is placed at the same ``level'' as $ n $).
\end{itemize}

\begin{prop}
    If $ r \in \mathbb{N}^+ $ and $ G \in \textrm{dom} ( \Phi ) $ are such that, for every $ v_0 , v_1 \in V ( G_n ) $, $ ( v_0 , v_1 ) \in E ( G_n ) $ if, and only if, there exists $ i \in \{ 0 , 1 \} $ such that $ l_r ( v_{1 - i} ) + 1 = l_r ( v_i ) $ (two joined vertices must be placed at consecutive ``levels'') and $ \left\lfloor \displaystyle\frac{v_i - m_r ( v_i )}{r} \right\rfloor = v_{1 - i} - m_r ( v_{1 - i} ) $, then $ ( \Phi ( G ) )_n = 2^{n - \lfloor n / r \rfloor} $, for every $ n \in \mathbb{N} $.
\end{prop}

\underline{Proof}

Let $ n \in \mathbb{N} $.

The last row of $ \textrm{adj} ( G_{n + 1} ) $ has exactly one nonzero entry, viz. $ \left( n + 2 , \left\lfloor \dfrac{n}{r} \right\rfloor + 1 \right) $.

Therefore, $ ( \Phi ( G ) )_n = \displaystyle\sum_{k = 1}^{n + 1} \left( \dfrac{( \textrm{adj} ( G_{n + 1} ) )_{n + 2 , k}}{2^{k - n - 1}} \right) = \dfrac{1}{2^{\lfloor n / r \rfloor - n}} = 2^{n - \lfloor n / r \rfloor} $. \hfill $ \square $

\begin{center}
    \begin{tabular}{|c|c|c|c|}
    \hline
    & $ G_6 $ & $ \textrm{adj} ( G_6 ) $ & \\
    \hline
    $ r = 1 $ & \raisebox{-40pt}{\begin{tikzpicture}[scale=0.4]
  \node at (0,6) (A) {$ 1 $};
  \node at (0,5) (B) {$ 2 $};
  \node at (0,4) (C) {$ 3 $};
  \node at (0,3) (D) {$ 4 $};
  \node at (0,2) (E) {$ 5 $};
  \node at (0,1) (F) {$ 6 $};
  \node at (0,0) (G) {$ 7 $};
  \path[]
        (A) edge node {} (B)
		(B) edge node {} (C)
        (C) edge node {} (D)
        (D) edge node {} (E)
        (E) edge node {} (F)
        (F) edge node {} (G);
    \end{tikzpicture}} & $ \left[ \begin{array}{ccccccc} * & * & * & * & * & * & * \\ 1 & * & * & * & * & * & * \\ 0 & 1 & * & * & * & * & * \\ 0 & 0 & 1 & * & * & * & * \\ 0 & 0 & 0 & 1 & * & * & * \\ 0 & 0 & 0 & 0 & 1 & * & * \\ 0 & 0 & 0 & 0 & 0 & 1 & * \end{array} \right] $ & $ \begin{array}{l} ( \Phi ( G ) )_0 = 1 \\ ( \Phi ( G ) )_1 = 1 \\ ( \Phi ( G ) )_2 = 1 \\ ( \Phi ( G ) )_3 = 1 \\ ( \Phi ( G ) )_4 = 1 \\ ( \Phi ( G ) )_5 = 1 \end{array} $ \\
    \hline
    $ r = 2 $ & \raisebox{-32pt}{\begin{tikzpicture}[]
  \node at (0,2) (A) {$ 1 $};
  \node at (0,1) (B) {$ 2 $};
  \node at (2,1) (C) {$ 3 $};
  \node at (0,0) (D) {$ 4 $};
  \node at (1,0) (E) {$ 5 $};
  \node at (2,0) (F) {$ 6 $};
  \node at (3,0) (G) {$ 7 $};
  \path[]
        (A) edge node {} (B)
		edge node {} (C)
		(B) edge node {} (D)
        edge node {} (E)
        (C) edge node {} (F)
        edge node {} (G);
    \end{tikzpicture}} & $ \left[ \begin{array}{ccccccc} * & * & * & * & * & * & * \\ 1 & * & * & * & * & * & * \\ 1 & 0 & * & * & * & * & * \\ 0 & 1 & 0 & * & * & * & * \\ 0 & 1 & 0 & 0 & * & * & * \\ 0 & 0 & 1 & 0 & 0 & * & * \\ 0 & 0 & 1 & 0 & 0 & 0 & * \end{array} \right] $ & $ \begin{array}{l} ( \Phi ( G ) )_0 = 1 \\ ( \Phi ( G ) )_1 = 2 \\ ( \Phi ( G ) )_2 = 2 \\ ( \Phi ( G ) )_3 = 4 \\ ( \Phi ( G ) )_4 = 4 \\ ( \Phi ( G ) )_5 = 8 \end{array} $ \\
    \hline
    $ r = 3 $ & \raisebox{-32pt}{\begin{tikzpicture}[]
  \node at (0,2) (A) {$ 1 $};
  \node at (0,1) (B) {$ 2 $};
  \node at (1,1) (C) {$ 3 $};
  \node at (2,1) (D) {$ 4 $};
  \node at (0,0) (E) {$ 5 $};
  \node at (1,0) (F) {$ 6 $};
  \node at (2,0) (G) {$ 7 $};
  \path[]
        (A) edge node {} (B)
		edge node {} (C)
		edge node {} (D)
		(B) edge node {} (E)
        edge node {} (F)
        edge node {} (G);
    \end{tikzpicture}} & $ \left[ \begin{array}{ccccccc} * & * & * & * & * & * & * \\ 1 & * & * & * & * & * & * \\ 1 & 0 & * & * & * & * & * \\ 1 & 0 & 0 & * & * & * & * \\ 0 & 1 & 0 & 0 & * & * & * \\ 0 & 1 & 0 & 0 & 0 & * & * \\ 0 & 1 & 0 & 0 & 0 & 0 & * \end{array} \right] $ & $ \begin{array}{l} ( \Phi ( G ) )_0 = 1 \\ ( \Phi ( G ) )_1 = 2 \\ ( \Phi ( G ) )_2 = 4 \\ ( \Phi ( G ) )_3 = 4 \\ ( \Phi ( G ) )_4 = 8 \\ ( \Phi ( G ) )_5 = 16 \end{array} $ \\
    \hline
    \end{tabular}
\end{center}

The next example shows the integer sequence that is associated to the Cayley graph of the free group on two generators.

\begin{prop}
    If $ G \in \textrm{dom} ( \Phi ) $ is such that $ E ( G_n ) = \{ ( 1 , 2 ) , ( 2 , 1 ) \} \cup \left\lbrace ( i , j ) \in V ( G_n )^2 \ | \ \left\lfloor \dfrac{i}{3} \right\rfloor = j \ \textrm{or} \ \left\lfloor \dfrac{j}{3} \right\rfloor = i \right\rbrace $, for every $ n \in \mathbb{N}^+ $, then $ ( \Phi ( G ) )_n = \left\lbrace \begin{array}{ll} 1 & \textrm{if} \ n = 0 \\ 2^{n - \lfloor ( n - 1 ) / 3 \rfloor} & \textrm{if} \ n > 0 \end{array} \right. $, for every $ n \in \mathbb{N} $.
\end{prop}

\underline{Proof}

It is clear that $ ( \Phi ( G ) )_0 = 1 $.

Let $ n \in \mathbb{N}^+ $.

The last row of $ \textrm{adj} ( G_{n + 1} ) $ has exactly one nonzero entry, viz. $ \left( n + 2 , \left\lfloor \dfrac{n - 1}{3} \right\rfloor + 1 \right) $.

Therefore, $ ( \Phi ( G ) )_n = \displaystyle\sum_{k = 1}^{n + 1} \left( \dfrac{( \textrm{adj} ( G_{n + 1} ) )_{n + 2 , k}}{2^{k - n - 1}} \right) = \dfrac{1}{2^{\lfloor ( n - 1 ) / 3 \rfloor - n}} = 2^{n - \lfloor ( n - 1 ) / 3 \rfloor} $. \hfill $ \square $

\begin{center}
    \begin{tabular}{|c|c|c|}
    \hline
    $ G_7 $ & $ \textrm{adj} ( G_7 ) $ & \\
    \hline
    \raisebox{-60pt}{\begin{tikzpicture}[scale=0.75]
  \node at (2,1) (A) {$ 1 $};
  \node at (2,3) (B) {$ 2 $};
  \node at (4,1) (C) {$ 3 $};
  \node at (2,-1) (D) {$ 4 $};
  \node at (0,1) (E) {$ 5 $};
  \node at (1,3) (F) {$ 6 $};
  \node at (2,4) (G) {$ 7 $};
  \node at (3,3) (H) {$ 8 $};
  \path[]
        (A) edge node {} (B)
		edge node {} (C)
		edge node {} (D)
		edge node {} (E)
		(B) edge node {} (F)
		edge node {} (G)
		edge node {} (H);
    \end{tikzpicture}} & $ \left[ \begin{array}{cccccccc} * & * & * & * & * & * & * & * \\ 1 & * & * & * & * & * & * & * \\ 1 & 0 & * & * & * & * & * & * \\ 1 & 0 & 0 & * & * & * & * & * \\ 1 & 0 & 0 & 0 & * & * & * & * \\ 0 & 1 & 0 & 0 & 0 & * & * & * \\ 0 & 1 & 0 & 0 & 0 & 0 & * & * \\ 0 & 1 & 0 & 0 & 0 & 0 & 0 & * \end{array} \right] $ & $ \begin{array}{l} ( \Phi ( G ) )_0 = 1 \\ ( \Phi ( G ) )_1 = 2 \\ ( \Phi ( G ) )_2 = 4 \\ ( \Phi ( G ) )_3 = 8 \\ ( \Phi ( G ) )_4 = 8 \\ ( \Phi ( G ) )_5 = 16 \\ ( \Phi ( G ) )_6 = 32 \end{array} $ \\
    \hline
    \end{tabular}
\end{center}

\subsection{Some interesting examples} \label{SubsectionInteresting}

The following example shows the integer sequence that is associated to divisibility's Hasse diagram.

\begin{thm}
    Let $ G \in \textrm{dom} ( \Phi ) $ such that, for every $ n \in \mathbb{N} $ and $ v_0 , v_1 \in V ( G_n ) $, $ ( v_0 , v_1 ) \in E ( G_n ) $ if, and only if, the following conditions hold:
    \begin{itemize}
        \item $ v_0 \neq v_1 $,
        \item there exists $ i \in \{ 0 , 1 \} $ such that $ v_i | v_{1 - i} $ and there does not exist $ w \in \mathbb{N}^+ \backslash \{ v_i , v_{1 - i} \} $ such that $ v_i \ | \ w \ | \ v_{1 - i} $.
    \end{itemize}
    
    Let also $ m \in \mathbb{N} $, $ r , \alpha_1 , ... , \alpha_r \in \mathbb{N}^+ $ and $ p_1 , ... , p_r $ pairwise different prime numbers such that $ \displaystyle\prod_{k = 1}^r ( p_k^{\alpha_k} ) = m + 2 $. Then $ ( \Phi ( G ) )_m = 2^{m + 1} \cdot \left\lbrace \begin{array}{cl} 2^{- p_1^{\alpha_1 - 1}} & \textrm{if} \ r = 1 \\ \sum_{k = 1}^r ( 2^{- p_k^{\alpha_k}} ) & \textrm{if} \ r > 1 \end{array} \right. $.
\end{thm}

\underline{Proof}

The last row of $ \textrm{adj} ( G_{m + 1} ) $ has exactly one nonzero entry (viz. $ ( m + 2 , p_1^{\alpha_1 - 1} ) $) if $ r = 1 $ and exactly $ r $ nonzero entries (viz. $ ( m + 2 , p_1^{\alpha_1} ) , ... , ( m + 2 , p_r^{\alpha_r} ) $) otherwise.

Therefore, $ ( \Phi ( G ) )_m = \displaystyle\sum_{k = 1}^{m + 1} \left( \dfrac{( \textrm{adj} ( G_{m + 1} ) )_{m + 2 , k}}{2^{k - m - 1}} \right) = \left\lbrace \begin{array}{cl} 1 / 2^{p_1^{\alpha_1 - 1} - m - 1} & \textrm{if} \ r = 1 \\ \sum_{k = 1}^r ( 1 / 2^{p_k^{\alpha_k} - m - 1} ) & \textrm{if} \ r > 1 \end{array} \right. = $

$ 2^{m + 1} \cdot \left\lbrace \begin{array}{cl} 2^{- p_1^{\alpha_1 - 1}} & \textrm{if} \ r = 1 \\ \sum_{k = 1}^r ( 2^{- p_k^{\alpha_k}} ) & \textrm{if} \ r > 1 \end{array} \right. $. \hfill $ \square $

\begin{center}
\begin{tabular}{|c|c|c|}
    \hline
    $ G_6 $ & $ \textrm{adj} ( G_6 ) $ & \\
    \hline
    \raisebox{-30pt}{\begin{tikzpicture}[]
  \node at (0,0) (A) {$ 1 $};
  \node at (0,1) (B) {$ 2 $};
  \node at (1,1) (C) {$ 3 $};
  \node at (0,2) (D) {$ 4 $};
  \node at (2,1) (E) {$ 5 $};
  \node at (1,2) (F) {$ 6 $};
  \node at (2,0) (G) {$ 7 $};
  \path[]
        (A) edge node {} (B)
		edge node {} (C)
		edge node {} (E)
		edge node {} (G)
		(B) edge node {} (D)
		edge node {} (F)
		(C) edge node {} (F);
    \end{tikzpicture}} & $ \left[ \begin{array}{ccccccc} * & * & * & * & * & * & * \\ 1 & * & * & * & * & * & * \\ 1 & 0 & * & * & * & * & * \\ 0 & 1 & 0 & * & * & * & * \\ 1 & 0 & 0 & 0 & * & * & * \\ 0 & 1 & 1 & 0 & 0 & * & * \\ 1 & 0 & 0 & 0 & 0 & 0 & * \end{array} \right] $ & $ \begin{array}{l} ( \Phi ( G ) )_0 = 1 \\ ( \Phi ( G ) )_1 = 2 \\ ( \Phi ( G ) )_2 = 2 \\ ( \Phi ( G ) )_3 = 8 \\ ( \Phi ( G ) )_4 = 12 \\ ( \Phi ( G ) )_5 = 32 \end{array} $ \\
    \hline
\end{tabular}
\end{center}

The following Maple functions help to compute $ \Phi ( G ) $, with $ G $ being as in the previous proposition.

\begin{code}{}
r:=m->numelems(ifactors(m+2)[2]):
p:=(k,m)->ifactors(m+2)[2][k][1]:
a:=(k,m)->ifactors(m+2)[2][k][2]:
Hasse:=m->2^(m+1)*piecewise((r(m)=1,2^(-p(1,m)^(a(1,m)-1))),sum(2^(-p(k,m)^a(k,m)),k=1..r(m))):
\end{code}

For example, the following command returns $ 8 $.

\begin{code}{}
Hasse(6);
\end{code}

The following problem exposes a possible pattern in divisibility’s Hasse diagram.

\begin{problem}
    If $ G $ is as in the previous proposition, what is the cardinality of
    
    $ \left\lbrace n \in \mathbb{N} \ | \ \dfrac{( \Phi ( G ) )_{n + 6}}{2^6} , \dfrac{( \Phi ( G ) )_{n + 5}}{2^4} , \dfrac{( \Phi ( G ) )_{n + 4}}{2^2} , ( \Phi ( G ) )_{n + 3} \geq 2^n \right\rbrace $?
\end{problem}

The following Maple program (which requires to have loaded the previous program) computes the elements of the previous problem's set that are upper-bounded by a given $ r \in \mathbb{N} $.

\begin{code}{}
test:=proc(r)
	local n, F:
	F:={}:
	for n from 0 to r do
		if Hasse(n+6)>=2^(n+6) and Hasse(n+5)>=2^(n+4) and Hasse(n+4)>=2^(n+2) and Hasse(n+3)>=2^n then
			F:=`union`(F,{n}):
		fi:
	od:
	return F:
end proc:
\end{code}

For example, the following command returns $ \{ 3 , 15 , 23 , 63 , 95 , 143 , 159 , 191 , 231 , 263 , 303 , 351 , 375 , 423 , 455 , 495 \} $.

\begin{code}{}
test(500);
\end{code}

\begin{problem}
   Find, if possible, an explicit formula for the strictly increasing integer sequence whose image is the previous problem's set, assuming that it has affirmative answer.
\end{problem}

The next example shows the integer sequence that is associated to Collatz's graph.

\begin{thm}
    If $ G \in \textrm{dom} ( \Phi ) $ is such that, for every $ n \in \mathbb{N} $ and $ v_0 , v_1 \in V ( G_n ) $, $ ( v_0 , v_1 ) \in E ( G_n ) $ if, and only if, there exists $ i \in \{ 0 , 1 \} $ such that $ \left\lbrace \begin{array}{c} v_{1 - i} = 3 \cdot v_i + 1 \\ v_i \ \textrm{is odd} \end{array} \right. $ or $ \left\lbrace \begin{array}{c} v_{1 - i} = v_i / 2 \\ v_i \ \textrm{is even} \end{array} \right. $, then $ ( \Phi ( G ) )_n = \left\lbrace \begin{array}{ll} 0 & \textrm{if} \ n \ \textrm{is odd} \\ \sqrt{2^n} + \sqrt[3]{4^{n + 1}} & \textrm{if} \ n \in \{ 6 \cdot r + 2 \}_{r \in \mathbb{N}} \\ \sqrt{2^n} & \textrm{otherwise} \end{array} \right. $, for every $ n \in \mathbb{N} $.
\end{thm}

\underline{Proof}

Note that no odd number is of the form $ 3 \cdot k + 1 $, with $ k $ being an odd number.

Let $ n \in \mathbb{N} $.

Then:
\begin{itemize}
    \item if $ n + 2 $ is odd, then the last row of $ \textrm{adj} ( G_{n + 1} ) $ has no nonzero entries,
    \item if $ n + 2 $ is even and of the form $ 3 \cdot k + 1 $, with $ k $ being an odd number, then the last row of $ \textrm{adj} ( G_{n + 1} ) $ has exactly two nonzero entries, viz. $ \left( n + 2 , \dfrac{n + 2}{2} \right) $ and $ \left( n + 2 , \dfrac{n + 1}{3} \right) $,
    \item in any other case, the last row of $ \textrm{adj} ( G_{n + 1} ) $ has exactly one nonzero entry, viz. $ \left( n + 2 , \dfrac{n + 2}{2} \right) $.
\end{itemize}

And note that $ n + 2 $ is of the form $ 3 \cdot k + 1 $, with $ k $ being an odd number, if, and only if, $ n \in \{ 6 \cdot r + 2 \}_{r \in \mathbb{N}} $.

Therefore, $ ( \Phi ( G ) )_n = \displaystyle\sum_{k = 1}^{n + 1} \left( \dfrac{( \textrm{adj} ( G_{n + 1} ) )_{n + 2 , k}}{2^{k - n - 1}} \right) = $

$ \left\lbrace \begin{array}{ll} 0 & \textrm{if} \ n \ \textrm{is odd} \\ 1 / 2^{( n + 2 ) / 2 - n - 1} + 1 / 2^{( n + 1 ) / 3 - n - 1} & \textrm{if} \ n \in \{ 6 \cdot r + 2 \}_{r \in \mathbb{N}} \\ 1 / 2^{( n + 2 ) / 2 - n - 1} & \textrm{otherwise} \end{array} \right. = \left\lbrace \begin{array}{ll} 0 & \textrm{if} \ n \ \textrm{is odd} \\ \sqrt{2^n} + \sqrt[3]{4^{n + 1}} & \textrm{if} \ n \in \{ 6 \cdot r + 2 \}_{r \in \mathbb{N}} \\ \sqrt{2^n} & \textrm{otherwise} \end{array} \right. $ \hfill $ \square $

\begin{center}
    \begin{tabular}{|c|c|c|}
    \hline
    $ G_{11} $ & $ \textrm{adj} ( G_{11} ) $ & \\
    \hline
    \raisebox{-45pt}{\begin{tikzpicture}[]
  \node at (0,3) (1) {$ 1 $};
  \node at (1,3) (2) {$ 2 $};
  \node at (2,3) (3) {$ 3 $};
  \node at (0,2) (4) {$ 4 $};
  \node at (1,2) (5) {$ 5 $};
  \node at (2,2) (6) {$ 6 $};
  \node at (0,0) (7) {$ 7 $};
  \node at (0,1) (8) {$ 8 $};
  \node at (1,0) (9) {$ 9 $};
  \node at (1,1) (10) {$ 10 $};
  \node at (2,0) (11) {$ 11 $};
  \node at (2,1) (12) {$ 12 $};
  \path[]
        (1) edge node {} (2)
		edge node {} (4)
		(3) edge node {} (6)
        edge node {} (10)
        (4) edge node {} (2)
        edge node {} (8)
        (5) edge node {} (10)
        (6) edge node {} (12);
    \end{tikzpicture}} & $ \left[ \begin{array}{cccccccccccc} * & * & * & * & * & * & * & * & * & * & * & * \\ 1 & * & * & * & * & * & * & * & * & * & * & * \\ 0 & 0 & * & * & * & * & * & * & * & * & * & * \\ 1 & 1 & 0 & * & * & * & * & * & * & * & * & * \\ 0 & 0 & 0 & 0 & * & * & * & * & * & * & * & * \\ 0 & 0 & 1 & 0 & 0 & * & * & * & * & * & * & * \\ 0 & 0 & 0 & 0 & 0 & 0 & * & * & * & * & * & * \\ 0 & 0 & 0 & 1 & 0 & 0 & 0 & * & * & * & * & * \\ 0 & 0 & 0 & 0 & 0 & 0 & 0 & 0 & * & * & * & * \\ 0 & 0 & 1 & 0 & 1 & 0 & 0 & 0 & 0 & * & * & * \\ 0 & 0 & 0 & 0 & 0 & 0 & 0 & 0 & 0 & 0 & * & * \\ 0 & 0 & 0 & 0 & 0 & 1 & 0 & 0 & 0 & 0 & 0 & * \end{array} \right] $ & $ \begin{array}{l} ( \Phi ( G ) )_0 = 1 \\ ( \Phi ( G ) )_1 = 0 \\ ( \Phi ( G ) )_2 = 6 \\ ( \Phi ( G ) )_3 = 0 \\ ( \Phi ( G ) )_4 = 4 \\ ( \Phi ( G ) )_5 = 0 \\ ( \Phi ( G ) )_6 = 8 \\ ( \Phi ( G ) )_7 = 0 \\ ( \Phi ( G ) )_8 = 80 \\ ( \Phi ( G ) )_9 = 0 \\ ( \Phi ( G ) )_{10} = 32 \end{array} $ \\
    \hline
\end{tabular}
\end{center}

The following Maple function computes $ \Phi ( G ) $, with $ G $ being as in the previous proposition.

\begin{code}{}
Collatz:=n->piecewise((gcd(n,2)=1,0),(gcd(n-2,6)=6,2^(2/3*(n+1))+2^(n/2)),2^(n/2)):
\end{code}

For example, the following command returns $ 0 $.

\begin{code}{}
Collatz(11);
\end{code}

\begin{problem}
    If $ G $ is as in the previous proposition and $ n \in \mathbb{N} $ such that $ n \geq 3 $, is the the number of connected components of $ G_n $ equal to $ \left\lfloor \displaystyle\frac{n + 2}{2} \right\rfloor - \left\lfloor \displaystyle\frac{n + 3}{6} \right\rfloor + 1 $?
\end{problem}

\begin{problem}[Collatz's Problem]
    If $ G $ is as in the previous proposition, is $ G_\infty $ connected?
\end{problem}

Regarding to the next two problems, recall the problems that are posed in Subsubsection \ref{Little}.

\begin{problem}
    Let $ D $ be the digraph such that $ v ( D ) = \mathbb{N}^+ $ and, for every $ m , n \in \mathbb{N}^+ $, $ ( m , n ) \in E ( D ) $ if, and only if, $ m \cdot 2^{n + 2} + 1 $ is a prime number that does not divide $ \dfrac{( F_{n + 2} - 1)^m - 1}{F_{n + 2} - 2} $. Let also $ H $ be the graph that is obtained by removing the loops of the underlying graph of the digraph $ D $. Find a vertex-by-vertex increasing graph sequence $ G $ such that $ G_\infty = H $ and, if possible, an explicit formula for $ \Phi ( G ) $.
\end{problem}

\begin{center}
\begin{tabular}{|c|c|c|}
    \hline
    $ G_5 $ & $ \textrm{adj} ( G_5 ) $ & \\
    \hline
    \raisebox{-5pt}{\begin{tikzpicture}[scale=0.75]
  \node at (0,1) (1) {$ 1 $};
  \node at (1,1) (2) {$ 2 $};
  \node at (2,1) (3) {$ 3 $};
  \node at (2,0) (4) {$ 4 $};
  \node at (1,0) (5) {$ 5 $};
  \node at (0,0) (6) {$ 6 $};
  \path[]
        (1) edge node {} (2)
		edge node {} (6)
		(2) edge node {} (5);
    \end{tikzpicture}} & $ \left[ \begin{array}{cccccc} * & * & * & * & * & * \\ 1 & * & * & * & * & * \\ 0 & 0 & * & * & * & * \\ 0 & 0 & 0 & * & * & * \\ 0 & 1 & 0 & 0 & * & * \\ 1 & 0 & 0 & 0 & 0 & * \end{array} \right] $ & $ \begin{array}{l} ( \Phi ( G ) )_0 = 1 \\ ( \Phi ( G ) )_1 = 0 \\ ( \Phi ( G ) )_2 = 0 \\ ( \Phi ( G ) )_3 = 4 \\ ( \Phi ( G ) )_4 = 16 \end{array} $ \\
    \hline
\end{tabular}
\end{center}

\begin{problem}
    If $ G $ is as in the previous problem, is $ G_\infty $ connected?
\end{problem}

The next problem asks for the integer sequence that is associated to the Hasse diagram that is induced by the Boolean algebra $ \wp ( \mathbb{N} ) $.

\begin{problem}
    Find, if possible, an explicit formula for $ \Phi ( G ) $, with $ G \in \textrm{dom} ( \Phi ) $ being such that, for every $ n \in \mathbb{N}^+ $ and $ v_0 , v_1 \in V ( G_n ) $, $ ( v_0 , v_1 ) \in E ( G_n ) $ if, and only if, there exist $ S_0 , S_1 \subseteq \{ 0 , ... , n \} $ such that:
    \begin{itemize}
        \item there exists $ i \in \{ 0 , 1 \} $ such that $ S_i \subsetneqq S_{1 - i} $ and there does not exist $ T \subsetneqq S_{1 - i} \ \textrm{such that} \ S_i \subsetneqq T $,
        \item $ \displaystyle\sum_{s \in S_k} ( 2^s ) + 1 = v_k $, for every $ k \in \{ 0 , 1 \} \ \left( \textrm{recall that} \ \displaystyle\sum_{s \in \varnothing} ( 2^s ) = 0 \right) $.
    \end{itemize}
\end{problem}

\begin{center}
\begin{tabular}{|c|c|c|}
    \hline
    $ G_7 $ & $ \textrm{adj} ( G_7 ) $ & \\
    \hline
    \raisebox{-45pt}{\begin{tikzpicture}[]
  \node at (1,0) (1) {$ 1 $};
  \node at (0,1) (2) {$ 2 $};
  \node at (1,1) (3) {$ 3 $};
  \node at (0,2) (4) {$ 4 $};
  \node at (2,1) (5) {$ 5 $};
  \node at (1,2) (6) {$ 6 $};
  \node at (2,2) (7) {$ 7 $};
  \node at (1,3) (8) {$ 8 $};
  \path[]
        (1) edge node {} (2)
		edge node {} (3)
		edge node {} (5)
		(2) edge node {} (4)
        edge node {} (6)
        (3) edge node {} (4)
        edge node {} (7)
        (4) edge node {} (8)
        (5) edge node {} (6)
        edge node {} (7)
        (6) edge node {} (8)
        (7) edge node {} (8);
    \end{tikzpicture}} & $ \left[ \begin{array}{cccccccc} * & * & * & * & * & * & * & * \\ 1 & * & * & * & * & * & * & * \\ 1 & 0 & * & * & * & * & * & * \\ 0 & 1 & 1 & * & * & * & * & * \\ 1 & 0 & 0 & 0 & * & * & * & * \\ 0 & 1 & 0 & 0 & 1 & * & * & * \\ 0 & 0 & 1 & 0 & 1 & 0 & * & * \\ 0 & 0 & 0 & 1 & 0 & 1 & 1 & * \end{array} \right] $ & $ \begin{array}{l} ( \Phi ( G ) )_0 = 1 \\ ( \Phi ( G ) )_1 = 2 \\ ( \Phi ( G ) )_2 = 3 \\ ( \Phi ( G ) )_3 = 8 \\ ( \Phi ( G ) )_4 = 9 \\ ( \Phi ( G ) )_5 = 10 \\ ( \Phi ( G ) )_6 = 11 \end{array} $ \\
    \hline
\end{tabular}
\end{center}

Note that, if $ G $ is as in the previous problem and $ n \in \mathbb{N} $, then $ ( \Phi ( G ) )_{2^n - 1} = 2^{2^n - 1} $ and $ G_{2^n - 1} $ is isomorphic to the $ n $-hypercube graph.

Apparently, previous problem's sequence coincides with OEIS A253317 from the second term of the later.

The following problem asks for the integer sequence that is associated to the graph encoding the proofs of the propositional logic.

\begin{problem}
    Let $ \phi $ be an enumeration of all theorems of the propositional logic such that $ \textrm{dom} ( \phi ) = \mathbb{N}^+ $ and $ \mathcal{G} \circ \phi $ is increasing, for some Gödel numbering $ \mathcal{G} $. Let also $ G \in \textrm{dom} ( \Phi ) $ such that, for every $ n \in \mathbb{N}^+ $ and $ i , j \in V ( G_n ) $, $ ( i , j ) \in E ( G_n ) $ if, and only if, there exist theorems $ \alpha , \beta $ and $ k \in V ( G_n ) $ such that $ \{ \phi_i , \phi_j \} = \{ [ \alpha \rightarrow \beta ] , \beta \} $ and $ \phi_k = \alpha $. Find, if possible, an explicit formula for $ \Phi ( G ) $.
\end{problem}

Recall that, given a digraph $ D $ and an equivalence relation $ \sim $ on $ D $, $ {}^{D}/_\sim $ denotes the graph whose vertex set is $ {}^{V ( D )}/_\sim $ and whose edge set is $ \{ ( {}^{v}/_\sim , {}^{w}/_\sim ) \ | \ v , w \in V ( D ) \ \textrm{and} \ ( v , w ) \in E ( D ) \} $.

\begin{problem}
    Let $ G $ be as in the previous problem, $ \sim \ = \{ ( i , j ) \in ( \mathbb{N}^+ )^2 \ | \ \phi_k = [ \phi_i \leftrightarrow \phi_j ] , \ \textrm{for some} \ k \in \mathbb{N}^+ \} $ and $ H \in \textrm{dom} ( \Phi ) $ such that $ H_\infty $ is isomorphic to $ {}^{G_\infty}/_\sim $. Find, if possible, an explicit formula for $ \Phi ( H ) $.
\end{problem}

An answer for the following problem would supply interesting information about the group $ \textrm{Sym}_\infty $.

\begin{problem}
    Let $ \phi $ be the enumeration of all finite permutations in reversed colexicographic ordering such that $ \textrm{dom} ( \phi ) = \mathbb{N}^+ $ (so, for example, $ ( \phi_1 , ... , \phi_6 ) = ( [ 1 ] , [ 2, 1 ] , [ 1, 3, 2 ] , [ 3, 1, 2 ] , [ 2, 3, 1 ] , [ 3, 2, 1 ] ) $, see OEIS A055089). Let also $ G \in \textrm{dom} ( \Phi ) $ such that, for every $ n \in \mathbb{N}^+ $ and $ i , j \in V ( G_n ) $, $ ( i , j ) \in E ( G_n ) $ if, and only if, $ \phi_i $ and $ \phi_j $ differ by exactly one transposition. Find, if possible, an explicit formula for $ \Phi ( G ) $.
\end{problem}

\begin{center}
\begin{tabular}{|c|c|c|}
    \hline
    $ G_5 $ & $ \textrm{adj} ( G_5 ) $ & \\
    \hline
    \raisebox{-45pt}{\begin{tikzpicture}[scale=0.75]
  \node at (0,-2) (1) {$ 1 $};
  \node at (-1.732050,-1) (2) {$ 2 $};
  \node at (1.732050,-1) (3) {$ 3 $};
  \node at (1.732050,1) (4) {$ 4 $};
  \node at (-1.732050,1) (5) {$ 5 $};
  \node at (0,2) (6) {$ 6 $};
  \path[]
        (1) edge node {} (2)
		edge node {} (3)
		edge node {} (6)
		(2) edge node {} (4)
        edge node {} (5)
        (3) edge node {} (4)
        edge node {} (5)
        (4) edge node {} (6)
        (5) edge node {} (6);
    \end{tikzpicture}} & $ \left[ \begin{array}{cccccc} * & * & * & * & * & * \\ 1 & * & * & * & * & * \\ 1 & 0 & * & * & * & * \\ 0 & 1 & 1 & * & * & * \\ 0 & 1 & 1 & 0 & * & * \\ 1 & 0 & 0 & 1 & 1 & * \end{array} \right] $ & $ \begin{array}{l} ( \Phi ( G ) )_0 = 1 \\ ( \Phi ( G ) )_1 = 2 \\ ( \Phi ( G ) )_2 = 3 \\ ( \Phi ( G ) )_3 = 6 \\ ( \Phi ( G ) )_4 = 19 \end{array} $ \\
    \hline
\end{tabular}
\end{center}

Note that, if $ G $ is as in the previous problem and $ n \in \mathbb{N}^+ $, then $ G_{n ! - 1} $ is isomorphic to the $ n $-transposition graph, $ \textrm{adj} ( G_{n ! - 1} ) $ is persymmetric (i.e., symmetric with respect to both diagonals) and the graph given by the edges and vertices of the $ n $-permutohedron (i.e., the Cayley graph that is generated by the adjacent transpositions of $ \textrm{Sym}_n $) is a subgraph of $ G_{n ! - 1} $.

\section{Acknowledgments}

The author thanks Professor Matthias Baaz for his kind guidance during the realization of this article.

This research has been partly supported by FWF Austria, project numbers P31063 and P31955.

\bibliographystyle{unsrt}  


\end{document}